\documentstyle[12pt, epic, eepic]{article}
\makeatletter
\def\section{\@startsection {section}{1}{\z@}{-3.5ex plus -1ex minus 
 -.2ex}{2.3ex plus .2ex}{\normalsize\bf}}
\def\@maketitle{\newpage
 \null
 \vskip 2em \begin{center}
 {\large \@title \par} \vskip 1.5em {\normalsize \lineskip .5em
\begin{tabular}[t]{c}\@author 
 \end{tabular}\par} 
 \vskip 1em {\normalsize \@date} \end{center}
 \par
 \vskip 1.5em} 
\makeatother

\title{{\bf Orbifold subfactors from Hecke algebras II}\\
---Quantum doubles and braiding---}
\author{{\sc David E. Evans}\\
Department of Mathematics\\
University of Wales, Swansea\\
Singleton Park, Swansea SA2 8PP, Wales, U.K.\\
\vphantom{X}\\
{\sc Yasuyuki Kawahigashi}\\
Department of Mathematical Sciences\\
University of Tokyo, Komaba, Tokyo, 153, JAPAN\\
e-mail: {\tt yasuyuki@ms.u-tokyo.ac.jp}}
\date{}

\begin{document}
\maketitle

\input amssym.def
\newsymbol\varnothing 203F
\newsymbol\rtimes 226F
\def\emptyset{\varnothing}

\def\qed{{\unskip\nobreak\hfil\penalty50
\hskip2em\hbox{}\nobreak\hfil\rm Q.E.D.
\parfillskip=0pt \finalhyphendemerits=0\par}\medskip}
\def\proof{{\noindent{\it Proof:\quad}}}

\def\tb{\mathbin{\bar\otimes}}
\def\equi{\sim}
\def\isom{\cong}
\def\ti{\tilde}
\def\lan{\langle}
\def\ran{\rangle}
\def\bd{\partial}

\def\Ad{{\hbox{Ad}}}
\def\ad{{\hbox{ad}}}
\def\Aut{{\hbox{Aut}}}
\def\dim{{\hbox{dim}}}
\def\Dom{{\hbox{Dom}}}
\def\End{{\hbox{End}}}
\def\Exp{{\hbox{Exp}}}
\def\gr{{\rm gr}}
\def\Hom{{\hbox{Hom}}}
\def\Int{{\hbox{Int}}}
\def\mod{{\hbox{mod}}}
\def\Orb{{\hbox{Orb}}}
\def\Out{{\hbox{Out}}}
\def\Path{{\hbox{Path}}}
\def\Proj{{\hbox{Proj}}}
\def\Sp{{\hbox{Sp}}}
\def\sp{{\hbox{sp}}}
\def\Str{{\hbox{String}}}
\def\String{{\hbox{String}}}
\def\supp{{\hbox{supp}}}
\def\tr{{\hbox{tr}}}
\def\Tr{{\hbox{Tr}}}
\def\Tube{{\hbox{Tube}\;}}

\def\C{{\bf C}}
\def\N{{\bf N}}
\def\Q{{\bf Q}}
\def\R{{\bf R}}
\def\Z{{\bf Z}}
\def\A{{\cal A}}
\def\E{{\cal E}}
\def\G{{\cal G}}
\def\H{{\cal H}}
\def\K{{\cal K}}
\def\T{{\cal T}}
\def\l{{\cal L}}
\def\M{{\cal M}}
\def\n{{\cal N}}
\def\p{{\cal P}}
\def\q{{\cal Q}}
\def\r{{\cal R}}
\def\ri{{\cal R_{0,1}}}
\def\U{{\cal U}}
\def\V{{\cal V}}
\def\z{{\cal Z}}

\def\a{\alpha}
\def\be{\beta}
\def\de{\delta}
\def\e{\varepsilon}
\def\epsilon{\varepsilon}
\def\ga{\gamma}
\def\Ga{\Gamma}
\def\la{\lambda}
\def\La{\Lambda}
\def\phi{\varphi}
\def\si{\sigma}
\def\th{\theta}
\def\t{\tau}
\def\om{\omega}
\def\Om{\Omega}

\newtheorem{theorem}{Theorem}[section]
\newtheorem{lemma}[theorem]{Lemma}
\newtheorem{corollary}[theorem]{Corollary}
\newtheorem{definition}[theorem]{Definition}
\newtheorem{assumption}[theorem]{Assumption}
\newtheorem{proposition}[theorem]{Proposition}
\newtheorem{remark}[theorem]{Remark}
\newtheorem{example}[theorem]{Example}

\begin{abstract}
A. Ocneanu has observed that
a mysterious orbifold phenomenon occurs in
the system of the $M_\infty$-$M_\infty$ bimodules of the
asymptotic inclusion, a subfactor analogue of the
quantum double, of the
Jones subfactor of type $A_{2n+1}$.

We show that this is a general phenomenon and identify some of his
orbifolds with the ones in our sense as subfactors given
as simultaneous fixed point algebras
by working on the Hecke algebra subfactors of type $A$ of Wenzl.
That is, we work on their asymptotic inclusions and
show that the $M_\infty$-$M_\infty$ bimodules are described by certain
orbifolds (with ghosts) for $SU(3)_{3k}$.
We actually compute several examples of the
(dual) principal graphs of the asymptotic inclusions.

As a corollary of the identification of Ocneanu's orbifolds with
ours, we show that a non-degenerate braiding exists on
the even vertices of $D_{2n}$, $n>2$.
\end{abstract}

\section{Introduction}
\label{intro}

In the theory of subfactors initiated by V. F. R. Jones in \cite{J},
Ocneanu's paragroup theory \cite{O1} is fundamental
in  descriptions of the combinatorial
structures arising from subfactors.  Ocneanu's construction
of the asymptotic
inclusions, introduced in \cite{O1}, has recently caught much
attention as a subfactor analogue of the quantum double
construction of Drinfel$'$d in \cite{D}.
(See \cite{E2}, \cite{EK3}, \cite{G2}, \cite{K2}, \cite{K3},
\cite{LR}, \cite{Masuda}, \cite{S2} on asymptotic inclusions.)

As noted by Ocneanu, if we start with a subfactor
$N\subset M=N\rtimes G$, where $N$ is a hyperfinite II$_1$
factor with a free action of a finite group $G$ on $N$,
then the resulting asymptotic inclusion $M\vee (M'\cap M_\infty)
\subset M_\infty$ gives the tensor category of the quantum
double of $G$ as that of the $M_\infty$-$M_\infty$ bimodules.
(See \cite[Section 12.8]{EK4} and \cite{KMY} for example.)
Since a paragroup, arising from
a subfactor, is a certain ``quantization'' of an ordinary
group, Ocneanu's construction of the asymptotic inclusion 
can be regarded as a subfactor analogue of the quantum double
construction.  (See \cite[Sections 12.6, 12.7, 13.5]{EK4} for
its relation to topological quantum field theory and rational
conformal field theory.)

The dual principal graphs of the asymptotic inclusions
are hard to compute, in general, while their principal graphs
are easy to compute, as in \cite{EK3}, \cite[Section 12.6]{EK4},
\cite{O2}, \cite{O3}, \cite{O4}, as long as we know the fusion rule
of the $M$-$M$ bimodules of the original subfactor $N\subset M$.
From the above viewpoint of the quantum double, it is the dual
principal graph, or the system of its even vertices, strictly
speaking, that is more important of the two graphs.
(In some sense, the principal graph represents just a double
``without quantum''.  See \cite[Section 12.6]{EK4}.)
So it would be interesting to have concrete descriptions of the
dual principal graphs (or their even vertices) of the asymptotic
inclusions of concrete examples of subfactors, other than the
ones of the form $N\subset M=N\rtimes G$ arising from genuine
groups.  Other ``easy'' examples
of subfactors of finite depth arising from actions of finite
groups contain
subgroup-group subfactors $N=R\rtimes H\subset M=R\rtimes G$
and Wassermann type subfactors $(\C\otimes R)^K\subset
(M_n(\C)\otimes R)^K$, where $R$ is the hyperfinite II$_1$ factor,
$H\subset G$ are finite groups acting freely on $R$, and
$K$ acts on $R=\bigotimes_k M_n(\C)$ as a product type action.
These, however, do not give anything new in the tensor categories
of their $M_\infty$-$M_\infty$ bimodules, because the
$M$-$M$ bimodules of these subfactors are given by the
tensor category of $\hat G$ and
then they give the same  $M_\infty$-$M_\infty$ bimodules as
the subfactor $N=R\rtimes H\subset M=R\rtimes G$, as seen from
\cite[Section 12.6]{EK4}.  In this sense, ``classical''
subfactors do not give interesting asymptotic inclusions.

The easiest subfactors among ``quantum'' subfactors
are the Jones subfactors $N\subset M$
of type $A_n$, as introduced in \cite{J}.
They are described as $N=\langle e_2, e_3, e_4,\dots\rangle$,
$M=\langle e_1, e_2, e_3, e_4,\dots\rangle$, where $\{e_j\}_{j\ge1}$
is a sequence of projections satisfying the following relations~:
\begin{eqnarray*}
e_j e_k&=& e_k e_j, \quad |j-k|\neq1,\\
e_j e_{j\pm1} e_j&=&(4\cos^2\frac{\pi}{n+1})^{-1} e_j.
\end{eqnarray*}
Then it is easy to see that the asymptotic inclusion
$M\vee (M'\cap M_\infty)\subset M_\infty$ is given as
\begin{eqnarray*}
M\vee (M'\cap M_\infty)&=&
\langle \dots, e_{-2}, e_{-1},e_1, e_2,\dots\rangle,\\
M_\infty&=&\langle \dots, e_{-2}, e_{-1},e_0, e_1, e_2,\dots\rangle,
\end{eqnarray*}
where $\{e_j\}_{j\in\Z}$ is a (double-sided) sequence of projections
satisfying the same relations as above.
The Jones indices of these subfactors were first computed by M. Choda
in \cite{Ch}.  It is quite non-trivial to describe the dual
principal graphs of these asymptotic inclusions, while the general
theory mentioned above gives the principal graphs easily from
the fusion rules.  Ocneanu has
announced a description of the $M_\infty$-$M_\infty$ bimodules
in \cite{O7} and at several other conferences.

If we apply the quantum double construction of Drinfel$'$d
to a Hopf algebra
already having an $R$-matrix, then the resulting Hopf algebra
is just a ``double'' of the original algebra and nothing interesting
happens in this procedure.  One might suspect we have something
similar and not interesting for these asymptotic inclusions of the
Jones subfactors, because the Jones subfactors correspond to
the quantum groups ${\cal U}_q (sl_2)$ in some sense.  This,
however, is not true.
We have a more subtle and interesting situation due to a certain
degeneracy condition
of the braiding in the sense of Ocneanu \cite{O7}.

The non-degeneracy condition of this kind was first introduced by
Reshetikhin--Turaev \cite{RT} in their construction of topological
invariants of 3-manifolds realizing the physical
prediction of Witten \cite{Wi1}.  From the topological viewpoint,
this condition is quite natural and it is this non-degeneracy
that leads the Turaev--Viro invariant \cite{TV} of a 3-manifold
being the square of the absolute value of the Reshetikhin--Turaev
invariant as in \cite{T}.  (See also \cite{O7} for an operator
algebraic account of this theorem of Turaev.)

Ocneanu has observed that a certain orbifold construction,
similar to the orbifold construction in our sense in \cite{EK1}
as simultaneous crossed products,
is invoked in the process of the asymptotic inclusions of the
Jones subfactors if the above non-degeneracy condition fails.
In this paper, we will generalize his construction to the case
of the Hecke algebra subfactors of Wenzl \cite{We} and show that
this is a general phenomenon in the following sense.  The asymptotic
inclusion produces a non-degenerate system of bimodules in
the sense of Ocneanu \cite{O7}.  From the viewpoint of
\cite{O7}, we can say that our orbifold construction \cite{EK1}
removes the degeneracy.  So if we apply the construction of the
asymptotic inclusion to a subfactor having a degenerate system
of bimodules, the orbifold construction is invoked automatically
in order to remove the degeneracy in the  procedure of making
the ``double''.  In this way, we get another series
of orbifold subfactors from Hecke akgebras of type $A$ as a
continuation of our work in \cite{EK1}.

The asymptotic inclusions of the Hecke algebra subfactors
are described naturally as follows.
The original subfactor of Wenzl is described as
$N=\langle g_2, g_3, g_4,\dots\rangle$,
$M=\langle g_1, g_2, g_3, g_4,\dots\rangle$, where $\{g_j\}_{j\ge1}$
is a sequence of the Hecke generators satisfying the relations
of the Hecke algebras of type $A$ as in \cite{We}.
The series of the commuting squares giving this subfactor
is not canonical in the sense of Popa, because this series
has a period larger than two.  Still, one can identify the
asymptotic inclusion of this subfactor as follows, which
is similar to the above description of the asymptotic
inclusions of the Jones subfactors~:
\begin{eqnarray*}
M\vee (M'\cap M_\infty)&=&
\langle \dots, g_{-2}, g_{-1},g_1, g_2,\dots\rangle,\\
M_\infty&=&\langle \dots, g_{-2}, g_{-1},g_0, g_1, g_2,\dots\rangle,
\end{eqnarray*}
where $\{g_j\}_{j\in\Z}$ is a (double-sided) sequence of the
Hecke generators satisfying the same relations.
This subfactor was first constructed with a double sided
sequence of the Hecke generators by Erlijman \cite{E}.
Later it was identified with the asymptotic inclusion
of the Hecke algebra subfactor of Wenzl by 
Goto \cite{G} and Erlijman \cite{E2} independently.  (Goto's
proof works in a quite general setting, while Erlijman directly
works on the Hecke algebras.)  So the asymptotic inclusions
of the Hecke algebra subfactors have natural constructions
in terms of generators and commuting squares parallel to 
the case of the Jones subfactors of type $A_n$.

In Section \ref{braid-tube}, we explain Ocneanu's
basic properties of braiding on a system of
bimodules in his sense \cite{O7} and its relation to his tube
algebras.  We continue the study of the tube algebras for
the Hecke algebra subfactors of Wenzl \cite{We} in
Section \ref{sect-braid}.  This gives the basic properties of the
tube algebra and enables us to use Ocneanu's general machinery
on asymptotic inclusions and tube algebras in \cite{O6}
and \cite{EK3}.
The dual principal graphs of the asymptotic inclusions
of the Jones subfactors of type $A_n$ are described in
Section \ref{sect-su2}.  This covers the case announced by
Ocneanu in \cite{O7}.  These for the Hecke algebra
subfactors with indices converging to 9 are dealt with in
Section \ref{sect-su3}.  This Section describes our main 
results.  In Section \ref{sect-orbif},
we study a relation between the orbifold phenomena Ocneanu
has observed and the orbifold construction in our sense \cite{EK1}
for the $SU(2)_{2k}$ case.
In the last Section \ref{sect-orbif-braid}, we study the
orbifold construction with braiding in our setting and get a
non-degenerate braiding on the even vertices of $D_{2n}$, $n>2$.

This work was done at University of Wales, Swansea
while the second author visited there on the joint research
program of the Royal Society and the Japan Society for Promotion
of Sciences.  The second author thanks for the financial
supports of these Societies.
The second author also acknowledges financial support
from the Inamori Foundation during this research. 
We thank Dr. Maxim Nazarov for his kind explanation of the
Littlewood--Richardson rule.

\section{Braiding and a tube algebra --- non-degenerate case ---}
\label{braid-tube}

We start with a finite braided system of bimodules
${\cal M}=\{x_i\}_{i\in I}$ in the sense of Ocneanu \cite{O7}.
(The original references for Ocneanu's theory used in this
paper are \cite{O1}, \cite{O2}, \cite{O3}, \cite{O4}, \cite{O5},
\cite{O6}, \cite{O7}.  See also \cite{EK2}, \cite{EK3} and
\cite[Chapter 12]{EK4}.)

An important example of such a system is
obtained from the WZW-models $SU(n)_k$ with Ocneanu's
surface bimodule construction as in \cite{O4}, \cite{O5}, \cite{O6}.
(See also \cite{EK3} or \cite[Chapter 12]{EK4}.)

We may have such a system from a subfactor $N\subset M$ with finite
index and finite depth.  Note that even when we have an abstract system of
bimodules, we can realize the system as a system of
bimodules arising from a single hyperfinite (possibly reducible)
subfactor $N\subset M$ of type II$_1$ finite index
and finite depth.  This is possible
by a minor variation of the construction in \cite{BG}.  That is, 
instead of choosing a primary field $\Phi$ in page 281 of \cite{BG},
we choose $\oplus_{i\in I} x_i$ as the generator to construct a paragroup.
In this way, we get a
(possibly reducible) subfactor $N\subset M$ for which the
system of the $M$-$M$ bimodules arising from the subfactor
is given by
${\cal M}=\{x_i\}_{i\in I}$.  (We have learnt this construction from
S. Yamagami.  See also \cite[Section 4]{EK2} or \cite[Section 12.5]{EK4}.)
So we may and do assume that our system ${\cal M}$ arises from a hyperfinite
subfactor $N\subset M$ of type II$_1$ with finite index and finite depth.

Define the {\sl global index}
$[{\cal M}]$ of the system by
$[{\cal M}]=\sum_{x\in {\cal M}} [x]$, where $[x]$ denotes the Jones
index of the bimodule $x$.
This is also the global index of the subfactor $N\subset M$ in the
sense of Ocneanu \cite{O1}.
We would like to study the system of the $M_\infty$-$M_\infty$
bimodules corresponding to the asymptotic inclusion
$M\vee (M'\cap M_\infty)\subset M_\infty$ arising from this
subfactor $N\subset M$ in the sense of Ocneanu \cite{O1}.
The construction of this system of the
$M_\infty$-$M_\infty$ bimodules can be regarded as
a subfactor analogue of the quantum double construction
of Drinfel$'$d in \cite{D}.  (This analogy has been noted
by Ocneanu.  See \cite[Chapter 12]{EK4} for the basic theory
of asymptotic inclusions and this analogy.)

In order to study this system, we work on Ocneanu's
tube algebra $\Tube {\cal M}$ and study the center of the tube algebra in the
sense of \cite{O4}, \cite{O7}.  (See also \cite[Chapter 12]{EK4} for
tube algebras.)

Recall that an element $x$ in the braided system
${\cal M}$ is called {\sl degenerate} in the sense of Ocneanu \cite{O7}
if it satisfies the identity in Figure \ref{degen}, where
the dashed circle
denotes the summation over all the labels $x\in{\cal M}$ with
coefficient $[x]^{1/2}/[{\cal M}]$ as in Figure \ref{killing}.
Such a dashed
ring is called a {\sl killing ring} in Ocneanu's terminology.

\unitlength 0.48mm
\thinlines
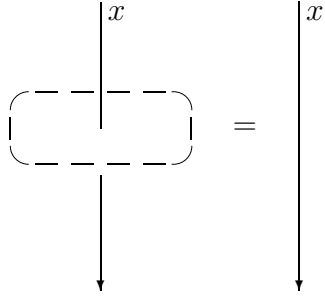
\begin{figure}[tb]
\begin{center}
\begin{picture}(135,100)
\put(25,40){\oval(10,10)[lb]}
\put(25,50){\oval(10,10)[lt]}
\put(65,40){\oval(10,10)[rb]}
\put(65,50){\oval(10,10)[rt]}
\put(20,42){\line(0,1){6}}
\put(70,42){\line(0,1){6}}
\multiput(27,35)(10,0){4}{\line(1,0){6}}
\multiput(27,55)(10,0){4}{\line(1,0){6}}
\put(45,32){\vector(0,-1){32}}
\put(45,45){\line(0,1){35}}
\put(85,45){\makebox(0,0){$=$}}
\put(100,80){\vector(0,-1){80}}
\put(47,75){\makebox(0,0)[lb]{$x$}}
\put(102,75){\makebox(0,0)[lb]{$x$}}
\end{picture}
\end{center}
\caption{A degenerate element $x$}
\label{degen}
\end{figure}

\unitlength 1.6mm
\thinlines
\begin{figure}[tb]
\begin{center}
\begin{picture}(50,20)
\put(40,15){\vector(-1,0){0}}
\put(40,10){\circle{10}}
\put(10,10){\arc{10}{0}{0.31416}}
\put(10,10){\arc{10}{0.62832}{0.94248}}
\put(10,10){\arc{10}{1.25664}{1.5708}}
\put(10,10){\arc{10}{1.88496}{2.19912}}
\put(10,10){\arc{10}{2.51328}{2.82744}}
\put(10,10){\arc{10}{3.1416}{3.45576}}
\put(10,10){\arc{10}{3.76992}{4.08408}}
\put(10,10){\arc{10}{4.39824}{4.7124}}
\put(10,10){\arc{10}{5.02656}{5.34072}}
\put(10,10){\arc{10}{5.65488}{5.96904}}
\put(26,10){\makebox(0,0){$=\displaystyle\sum_{x\in{\cal M}}
\frac{[x]^{1/2}}{[{\cal M}]}$}}
\put(38,16){\makebox(0,0){$x$}}
\end{picture}
\end{center}
\caption{A killing ring}
\label{killing}
\end{figure}

In this section, we suppose that the braiding on $\cal M$
is non-degenerate in the
sense that $0$ is the only degenerate element.  (We remark that
$0$ is always degenerate by definition.)
We note that we can use a graphical expression as in
\cite{KL}, \cite{Wi2} for elements in the tube algebra $\Tube {\cal M}$
because we have a braided system of bimodules.
In the tube algebra, we define the
{\sl Ocneanu projection} $p_{a,b}\in \Tube{\cal M}$ for $a,b\in {\cal M}$
as in Figure \ref{O-proj}.  In this picture,
the left half is a coefficient represented diagramatically  as
in \cite{KL}, where the horizontal bar represents a fraction,
and the right half is an element in the
tube algebra $\Tube {\cal M}$ where the top and the bottom of the dashed
line are identified in the tube picture.  The dashed line
again denotes the killing ring.

\unitlength 0.7mm
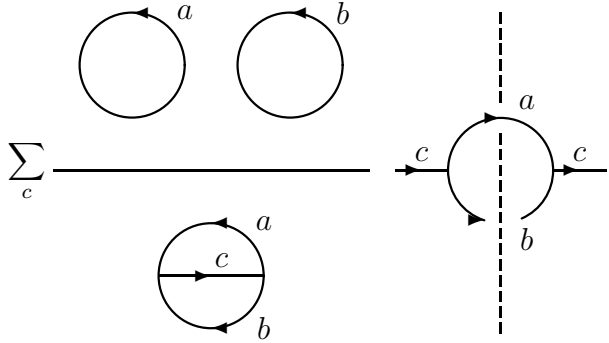
\begin{figure}[tb]
\begin{center}
\begin{picture}(120,80)
\thicklines
\put(20,60){\circle{20}}
\put(50,60){\circle{20}}
\put(35,20){\circle{20}}
\put(5,40){\line(1,0){60}}
\put(25,20){\line(1,0){20}}
\put(70,40){\line(1,0){10}}
\put(100,40){\line(1,0){10}}
\put(20,70){\vector(-1,0){0}}
\put(50,70){\vector(-1,0){0}}
\put(35,10){\vector(-1,0){0}}
\put(35,20){\vector(1,0){0}}
\put(35,30){\vector(-1,0){0}}
\put(75,40){\vector(1,0){0}}
\put(105,40){\vector(1,0){0}}
\put(90,50){\vector(1,0){0}}
\put(86.91,30.49){\vector(1,0){0}}
\put(90,40){\arc{20}{0}{1.157}}
\put(90,40){\arc{20}{1.885}{6.2832}}
\put(30,70){\makebox(0,0){$a$}}
\put(60,70){\makebox(0,0){$b$}}
\put(45,30){\makebox(0,0){$a$}}
\put(37,23){\makebox(0,0){$c$}}
\put(45,10){\makebox(0,0){$b$}}
\put(95,53){\makebox(0,0){$a$}}
\put(75,43){\makebox(0,0){$c$}}
\put(105,43){\makebox(0,0){$c$}}
\put(95,27){\makebox(0,0){$b$}}
\multiput(90,70)(0,-3){6}{\line(0,-1){2}}
\multiput(90,47)(0,-3){13}{\line(0,-1){2}}
\put(0,40){\makebox(0,0){$\displaystyle \sum_c$}}
\end{picture}
\end{center}
\caption{The Ocneanu projection $p_{a,b}$}
\label{O-proj}
\end{figure}

We recall from \cite[Section 12.3]{KL}
that we can perform the graphical operation
called a {\sl handle slide} against a killing ring without changing the
number or operator represented by the figure.
We give an
example of a handle slide in Figures \ref{handle1}, \ref{handle2}.
In this situation here, we assume that the link components on the
right hand side are killing rings.
(We remark that we have to regard a diagram of a link as
a {\sl framed} link diagram now.)
Note that
this handle slide is valid regardless of the non-degeneracy condition.
(See \cite[Section 12.3]{KL}.)
\thicklines
\unitlength 1mm
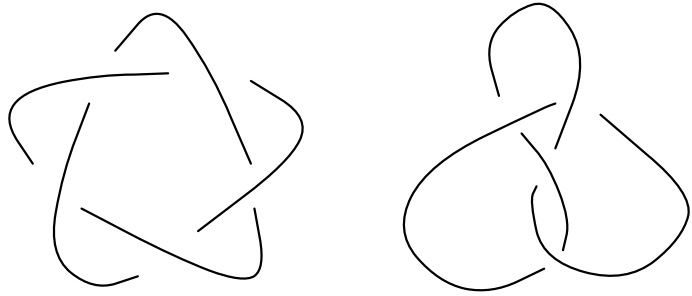
\begin{figure}[tb]
\begin{center}
\begin{picture}(120,60)
\spline(27,45)(33,52)(40,42)(45,30)
\spline(45,41)(53,36)(50,30)(38,21)
\spline(45.5,24)(47,16)(44,14)(35,17.5)(22.5,24)
\spline(30,15)(24,13)(18,18)(20,29)(23.5,38)
\spline(16,30)(12,36)(15,40)(25,42)(34,42)
\spline(78,39)(76,46)(80,50.5)(85,52)(90,44)(85.5,32)
\spline(83,27)(82,25)(84,17)(95,14)(102.5,20)(104,26)(91.5,36.5)
\spline(85.5,38)(84,37.5)(68.5,30)(64,21)(70,14)(77,12.5)(84,16)
\spline(86.5,18.5)(87.5,22.5)(84.5,30)(81,34)
\end{picture}
\end{center} 
\caption{Before a handle slide}
\label{handle1}
\end{figure}

\thicklines
\unitlength 1mm
\begin{figure}[tb]
\begin{center}
\begin{picture}(120,60)
\spline(27,45)(33,52)(40,42)(45,30)
\spline(45.5,24)(47,16)(44,14)(35,17.5)(22.5,24)
\spline(30,15)(24,13)(18,18)(20,29)(23.5,38)
\spline(16,30)(12,36)(15,40)(25,42)(34,42)
\spline(78,39)(76,46)(80,50.5)(85,52)(90,44)(85.5,32)
\spline(83,27)(82,25)(84,17)(95,14)(102.5,20)(104,26)(91.5,36.5)
\spline(85.5,38)(84,37.5)(68.5,30)(64,21)(70,14)(77,12.5)(84,16)
\spline(86.5,18.5)(87.5,22.5)(84.5,30)(81,34)
\spline(45,41)(60,37)(69,35)(79,37)(85,40)
\spline(92,37.5)(92.5,36.5)(105,26)(103.5,20)(95,13)(84,16)(81,20)(81,28)
\spline(79,39)(77,46)(80,49.5)(85,50.5)(89,44)(84.5,34)
\spline(87.5,18)(88.5,22.5)(85.5,30)(81,35)
\spline(85,15)(77,11.5)(70,11.5)(59,27)(51,28)(45,26)(38,21)
\end{picture}
\end{center} 
\caption{After a handle slide}
\label{handle2}
\end{figure}
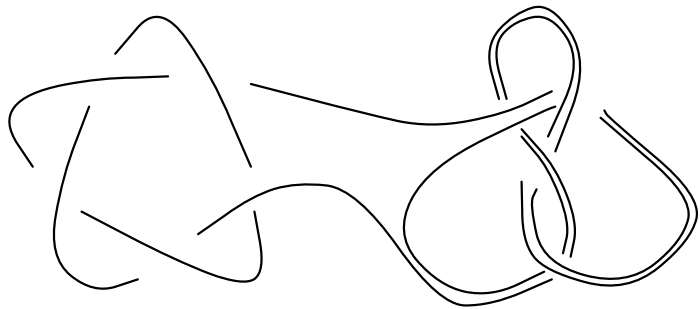

The following theorem is due to Ocneanu.  He presented this theorem and
a proof in his talk in the Taniguchi Symposium in Japan in July, 1993.
The proof here, except for the last paragraph, is his proof, which
we include for the sake of completeness.

\begin{theorem}
\label{ocn}
The above element $p_{a,b}$ gives a system of mutually orthogonal
minimal central projections in
the tube algebra with $\sum_{a,b\in{\cal M}} p_{a,b}=1$.
\end{theorem}

Before the proof, note that
this system of the minimal central
projections describes the system
of the $M_\infty$-$M_\infty$ bimodules and that these bimodules give the
even vertices of the dual principal graph of the asymptotic inclusion,
by Ocneanu's theorem in \cite{O6} (see \cite[Theorem 4.3]{EK3} or
\cite[Theorem 12.28]{EK4}).

\begin{proof}
First we prove that $p_{a,b}$'s give a system of mutually orthogonal
projections.  It is clear that each $p_{a,b}$ is self-adjoint, so we
will prove that $p_{a,b}p_{a',b'}=\delta_{a,a'}\delta_{b,b'}p_{a,b}$
first.  This proof is given as in Figure \ref{nondeg1},
where we compute $p_{a,b}p_{a',b'}$ graphically.  In the first
equality, we have used a handle slide.  In the third equality,
the non-degeneracy assumption implies that the label $d$ should be 0.
\unitlength 0.5mm
\begin{figure}[tb]
\begin{center}
\begin{picture}(260,320)
\thicklines
\put(50,60){\circle{20}}
\put(80,60){\circle{20}}
\put(65,20){\circle{20}}
\put(20,140){\circle{20}}
\put(50,140){\circle{20}}
\put(80,140){\circle{20}}
\put(110,140){\circle{20}}
\put(140,140){\circle{20}}
\put(35,100){\circle{20}}
\put(65,100){\circle{20}}
\put(95,100){\circle{20}}
\put(125,100){\circle{20}}
\put(20,220){\circle{20}}
\put(50,220){\circle{20}}
\put(35,180){\circle{20}}
\put(80,220){\circle{20}}
\put(110,220){\circle{20}}
\put(95,180){\circle{20}}
\put(20,300){\circle{20}}
\put(50,300){\circle{20}}
\put(35,260){\circle{20}}
\put(80,300){\circle{20}}
\put(110,300){\circle{20}}
\put(95,260){\circle{20}}
\put(35,40){\line(1,0){60}}
\put(55,20){\line(1,0){20}}
\put(100,40){\line(1,0){10}}
\put(130,40){\line(1,0){10}}
\put(5,120){\line(1,0){180}}
\put(25,100){\line(1,0){20}}
\put(55,100){\line(1,0){20}}
\put(85,100){\line(1,0){20}}
\put(115,100){\line(1,0){20}}
\put(190,120){\line(1,0){10}}
\put(240,120){\line(1,0){10}}
\put(5,200){\line(1,0){120}}
\put(25,180){\line(1,0){20}}
\put(85,180){\line(1,0){20}}
\put(130,200){\line(1,0){10}}
\put(160,200){\line(1,0){20}}
\put(200,200){\line(1,0){10}}
\put(5,280){\line(1,0){120}}
\put(25,260){\line(1,0){20}}
\put(85,260){\line(1,0){20}}
\put(130,280){\line(1,0){10}}
\put(160,280){\line(1,0){20}}
\put(200,280){\line(1,0){10}}
\put(50,70){\vector(-1,0){0}}
\put(80,70){\vector(-1,0){0}}
\put(20,150){\vector(-1,0){0}}
\put(50,150){\vector(-1,0){0}}
\put(80,150){\vector(-1,0){0}}
\put(110,150){\vector(-1,0){0}}
\put(140,150){\vector(-1,0){0}}
\put(20,230){\vector(-1,0){0}}
\put(50,230){\vector(-1,0){0}}
\put(80,230){\vector(-1,0){0}}
\put(110,230){\vector(-1,0){0}}
\put(20,310){\vector(-1,0){0}}
\put(50,310){\vector(-1,0){0}}
\put(80,310){\vector(-1,0){0}}
\put(110,310){\vector(-1,0){0}}
\put(65,10){\vector(-1,0){0}}
\put(65,20){\vector(1,0){0}}
\put(65,30){\vector(-1,0){0}}
\put(35,90){\vector(-1,0){0}}
\put(35,100){\vector(1,0){0}}
\put(35,110){\vector(-1,0){0}}
\put(65,90){\vector(-1,0){0}}
\put(65,100){\vector(1,0){0}}
\put(65,110){\vector(-1,0){0}}
\put(95,90){\vector(-1,0){0}}
\put(95,100){\vector(1,0){0}}
\put(95,110){\vector(-1,0){0}}
\put(125,90){\vector(-1,0){0}}
\put(125,100){\vector(1,0){0}}
\put(125,110){\vector(-1,0){0}}
\put(35,170){\vector(-1,0){0}}
\put(35,180){\vector(1,0){0}}
\put(35,190){\vector(-1,0){0}}
\put(95,170){\vector(-1,0){0}}
\put(95,180){\vector(1,0){0}}
\put(95,190){\vector(-1,0){0}}
\put(35,250){\vector(-1,0){0}}
\put(35,260){\vector(1,0){0}}
\put(35,270){\vector(-1,0){0}}
\put(95,250){\vector(-1,0){0}}
\put(95,260){\vector(1,0){0}}
\put(95,270){\vector(-1,0){0}}
\put(110,40){\vector(1,0){0}}
\put(140,40){\vector(1,0){0}}
\put(120,50){\vector(1,0){0}}
\put(116.91,30.49){\vector(1,0){0}}
\put(200,120){\vector(1,0){0}}
\put(250,120){\vector(1,0){0}}
\put(220,130){\vector(1,0){0}}
\put(220,110){\vector(1,0){0}}
\put(220,120){\vector(0,-1){0}}
\put(175,140){\vector(1,0){0}}
\put(185,140){\vector(0,-1){0}}
\put(170,150){\vector(-1,2){0}}
\put(170,130){\vector(-1,-2){0}}
\put(165,140){\vector(-1,2){0}}
\put(165,140){\vector(-1,-2){0}}
\put(165,140){\line(1,0){10}}
\put(165,140){\line(1,2){5}}
\put(165,140){\line(1,-2){5}}
\put(175,140){\line(-1,2){5}}
\put(175,140){\line(-1,-2){5}}
\put(140,200){\vector(1,0){0}}
\put(170,200){\vector(1,0){0}}
\put(210,200){\vector(1,0){0}}
\put(150,210){\vector(1,0){0}}
\put(190,210){\vector(1,0){0}}
\put(146.91,190.49){\vector(1,0){0}}
\put(190,190){\vector(1,0){0}}
\put(140,280){\vector(1,0){0}}
\put(170,280){\vector(1,0){0}}
\put(210,280){\vector(1,0){0}}
\put(150,290){\vector(1,0){0}}
\put(190,290){\vector(1,0){0}}
\put(146.91,270.49){\vector(1,0){0}}
\put(186.91,270.49){\vector(1,0){0}}
\put(210,130){\line(1,0){20}}
\put(213,110){\line(1,0){17}}
\put(220,110){\line(0,1){5}}
\put(220,120){\line(0,1){10}}
\put(120,40){\arc{20}{0}{1.157}}
\put(120,40){\arc{20}{1.885}{6.2832}}
\put(230,120){\arc{20}{0}{1.5708}}
\put(230,120){\arc{20}{4.7124}{6.2832}}
\put(175,140){\arc{20}{0}{1.5708}}
\put(175,140){\arc{20}{4.7124}{6.2832}}
\put(170,150){\line(1,0){5}}
\put(170,130){\line(1,0){5}}
\put(210,120){\arc{20}{1.885}{4.7124}}
\put(150,200){\arc{20}{0}{0.3}}
\put(150,200){\arc{20}{0.9}{1.2}}
\put(150,200){\arc{20}{1.885}{6.2832}}
\put(190,200){\arc{20}{0}{1.85}}
\put(190,200){\arc{20}{2.84}{6.2832}}
\put(150,280){\arc{20}{0}{1.157}}
\put(150,280){\arc{20}{1.885}{6.2832}}
\put(190,280){\arc{20}{0}{1.157}}
\put(190,280){\arc{20}{1.885}{6.2832}}
\multiput(120,70)(0,-3){6}{\line(0,-1){2}}
\multiput(120,47)(0,-3){13}{\line(0,-1){2}}
\multiput(210,150)(0,-3){6}{\line(0,-1){2}}
\multiput(210,127)(0,-3){13}{\line(0,-1){2}}
\multiput(150,230)(0,-3){6}{\line(0,-1){2}}
\multiput(150,207)(0,-3){13}{\line(0,-1){2}}
\multiput(150,310)(0,-3){6}{\line(0,-1){2}}
\multiput(150,287)(0,-3){13}{\line(0,-1){2}}
\multiput(190,310)(0,-3){6}{\line(0,-1){2}}
\multiput(190,287)(0,-3){13}{\line(0,-1){2}}
\put(60,70){\makebox(0,0){$a$}}
\put(90,70){\makebox(0,0){$b$}}
\put(75,30){\makebox(0,0){$a$}}
\put(67,23){\makebox(0,0){$c$}}
\put(75,10){\makebox(0,0){$b$}}
\put(125,53){\makebox(0,0){$a$}}
\put(105,44){\makebox(0,0){$c$}}
\put(135,44){\makebox(0,0){$c$}}
\put(125,27){\makebox(0,0){$b$}}
\put(30,150){\makebox(0,0){$a$}}
\put(60,150){\makebox(0,0){$b$}}
\put(92,150){\makebox(0,0){$a'$}}
\put(122,150){\makebox(0,0){$b'$}}
\put(150,150){\makebox(0,0){$d$}}
\put(45,110){\makebox(0,0){$a$}}
\put(37,103){\makebox(0,0){$c$}}
\put(45,90){\makebox(0,0){$b$}}
\put(77,110){\makebox(0,0){$a'$}}
\put(67,103){\makebox(0,0){$c$}}
\put(77,90){\makebox(0,0){$b'$}}
\put(105,110){\makebox(0,0){$b$}}
\put(97,104){\makebox(0,0){$b'$}}
\put(105,90){\makebox(0,0){$d$}}
\put(135,110){\makebox(0,0){$a$}}
\put(127,104){\makebox(0,0){$a'$}}
\put(135,90){\makebox(0,0){$d$}}
\put(30,230){\makebox(0,0){$a$}}
\put(60,230){\makebox(0,0){$b$}}
\put(92,230){\makebox(0,0){$a'$}}
\put(122,230){\makebox(0,0){$b'$}}
\put(30,310){\makebox(0,0){$a$}}
\put(60,310){\makebox(0,0){$b$}}
\put(92,310){\makebox(0,0){$a'$}}
\put(122,310){\makebox(0,0){$b'$}}
\put(45,190){\makebox(0,0){$a$}}
\put(37,183){\makebox(0,0){$c$}}
\put(45,170){\makebox(0,0){$b$}}
\put(106,190){\makebox(0,0){$a'$}}
\put(97,183){\makebox(0,0){$c$}}
\put(105,170){\makebox(0,0){$b'$}}
\put(45,270){\makebox(0,0){$a$}}
\put(37,263){\makebox(0,0){$c$}}
\put(45,250){\makebox(0,0){$b$}}
\put(106,270){\makebox(0,0){$a'$}}
\put(97,263){\makebox(0,0){$c$}}
\put(105,250){\makebox(0,0){$b'$}}
\put(155,213){\makebox(0,0){$a$}}
\put(135,204){\makebox(0,0){$c$}}
\put(170,204){\makebox(0,0){$c$}}
\put(155,186){\makebox(0,0){$b$}}
\put(195,213){\makebox(0,0){$a'$}}
\put(205,204){\makebox(0,0){$c$}}
\put(195,186){\makebox(0,0){$b'$}}
\put(155,293){\makebox(0,0){$a$}}
\put(135,284){\makebox(0,0){$c$}}
\put(170,284){\makebox(0,0){$c$}}
\put(155,267){\makebox(0,0){$b$}}
\put(195,293){\makebox(0,0){$a'$}}
\put(205,284){\makebox(0,0){$c$}}
\put(195,267){\makebox(0,0){$b'$}}
\put(162,145){\makebox(0,0){$a$}}
\put(162,135){\makebox(0,0){$b$}}
\put(177,145){\makebox(0,0){$a'$}}
\put(177,135){\makebox(0,0){$b'$}}
\put(170,143){\makebox(0,0){$c$}}
\put(189,140){\makebox(0,0){$d$}}
\put(195,124){\makebox(0,0){$c$}}
\put(245,124){\makebox(0,0){$c$}}
\put(202,133){\makebox(0,0){$a$}}
\put(202,107){\makebox(0,0){$b$}}
\put(238,133){\makebox(0,0){$a'$}}
\put(238,107){\makebox(0,0){$b'$}}
\put(224,125){\makebox(0,0){$d$}}
\multiput(159,207)(3.33,0){7}{\line(1,0){2}}
\multiput(158,193)(3.33,0){8}{\line(1,0){2}}
\put(155,200){\line(0,1){2}}
\put(155,204){\line(1,2){1}}
\put(155,197){\line(0,1){2}}
\put(155,196){\line(1,-2){1}}
\put(185,204){\line(-1,2){1}}
\put(185,200){\line(0,1){2}}
\put(185,197){\line(0,1){2}}
\put(185,196){\line(-1,-2){1}}
\put(219,117){\line(1,0){2}}
\put(222,118){\line(1,1){2}}
\put(217,118){\line(-1,1){2}}
\put(215,121){\line(0,1){2}}
\put(215,125){\line(1,1){2}}
\put(221,127){\line(1,0){2}}
\put(0,40){\makebox(0,0){$=\displaystyle
\delta_{aa'}\delta_{bb'}\sum_c$}}
\put(-6,120){\makebox(0,0){$=\displaystyle\sum_{c,d}$}}
\put(-6,200){\makebox(0,0){$=\displaystyle\sum_{c}$}}
\put(-6,280){\makebox(0,0){$\hphantom{=}\displaystyle\sum_c$}}
\end{picture}
\end{center}
\caption{Orthogonality of $p_{a,b}$}
\label{nondeg1}
\end{figure}
We next prove that each $p_{a,b}$ is central.  It is enough to prove
that $p_{a,b}/([a][b])^{1/2}$ commutes with any element in the tube
algebra.
This proof is given again graphically.  First we compute
the product of $p_{a,b}/([a][b])^{1/2}$ and a generic element in 
the tube algebra  as in Figure \ref{nondeg2}, where the top
and the bottom of the generic element labeled with $z$ are
identified again, where the top
and the bottom of the generic element labeled with $z$ are
identified again.
\unitlength 0.55mm
\begin{figure}[tb]
\begin{center}
\begin{picture}(160,320)
\thicklines
\put(20,20){\circle{20}}
\put(50,20){\circle{20}}
\put(20,100){\circle{20}}
\put(50,100){\circle{20}}
\put(35,140){\circle{20}}
\put(20,180){\circle{20}}
\put(20,260){\circle{20}}
\put(10,20){\vector(1,0){20}}
\put(40,20){\vector(1,0){20}}
\put(10,100){\vector(1,0){20}}
\put(40,100){\vector(1,0){20}}
\put(10,180){\vector(1,0){20}}
\put(10,260){\vector(1,0){20}}
\put(5,40){\line(1,0){60}}
\put(5,120){\line(1,0){60}}
\put(5,200){\line(1,0){30}}
\put(5,280){\line(1,0){30}}
\put(35,55){\makebox(0,0){$1$}}
\put(20,215){\makebox(0,0){$1$}}
\put(20,295){\makebox(0,0){$1$}}
\put(90,40){\vector(1,0){10}}
\put(110,40){\vector(1,0){10}}
\put(140,40){\vector(1,0){10}}
\put(70,120){\vector(1,0){10}}
\put(100,120){\vector(1,0){10}}
\put(130,120){\vector(1,0){10}}
\put(130,130){\vector(1,0){20}}
\put(40,200){\vector(1,0){10}}
\put(70,200){\vector(1,0){10}}
\put(70,210){\vector(1,0){20}}
\put(40,280){\vector(1,0){10}}
\put(70,280){\vector(1,0){10}}
\put(90,280){\vector(1,0){10}}
\put(100,290){\vector(1,0){10}}
\put(20,30){\vector(-1,0){0}}
\put(50,30){\vector(-1,0){0}}
\put(20,10){\vector(-1,0){0}}
\put(50,10){\vector(-1,0){0}}
\put(20,110){\vector(-1,0){0}}
\put(50,110){\vector(-1,0){0}}
\put(20,90){\vector(-1,0){0}}
\put(50,90){\vector(-1,0){0}}
\put(35,150){\vector(-1,0){0}}
\put(20,190){\vector(-1,0){0}}
\put(20,170){\vector(-1,0){0}}
\put(20,270){\vector(-1,0){0}}
\put(20,250){\vector(-1,0){0}}
\put(130,50){\vector(1,0){0}}
\put(126.91,30.49){\vector(1,0){0}}
\put(90,130){\vector(1,0){0}}
\put(86.91,110.49){\vector(1,0){0}}
\put(60,210){\vector(1,0){0}}
\put(56.91,190.49){\vector(1,0){0}}
\put(60,290){\vector(1,0){0}}
\put(56.91,270.49){\vector(1,0){0}}
\put(100,300){\vector(0,-1){40}}
\put(70,190){\vector(-1,0){0}}
\put(80,200){\vector(0,-1){0}}
\put(130,110){\vector(-1,0){0}}
\put(140,120){\vector(0,-1){0}}
\put(120,110){\vector(1,0){0}}
\put(120,130){\vector(1,0){0}}
\put(80,30){\vector(1,0){0}}
\put(80,50){\vector(1,0){0}}
\put(90,30){\vector(-1,0){0}}
\put(100,40){\vector(0,-1){0}}
\put(82.5,55){\vector(-1,0){0}}
\put(60,280){\arc{20}{0}{1.157}}
\put(60,280){\arc{20}{1.885}{6.2832}}
\put(90,120){\arc{20}{0}{1.157}}
\put(90,120){\arc{20}{1.885}{6.2832}}
\put(130,40){\arc{20}{0}{1.157}}
\put(130,40){\arc{20}{1.885}{6.2832}}
\put(60,200){\arc{20}{0}{0.75}}
\put(60,200){\arc{20}{1.15}{1.3}}
\put(60,200){\arc{20}{1.885}{6.2832}}
\put(70,200){\arc{20}{0}{3.9}}
\put(70,200){\arc{20}{4.3}{6.2832}}
\put(120,120){\arc{20}{0}{0.75}}
\put(120,120){\arc{20}{1.15}{6.2832}}
\put(130,120){\arc{20}{0}{3.9}}
\put(130,120){\arc{20}{4.3}{6.2832}}
\put(80,40){\arc{20}{0}{0.75}}
\put(80,40){\arc{20}{1.15}{6.2832}}
\put(90,40){\arc{20}{0}{3.9}}
\put(90,40){\arc{20}{4.3}{6.2832}}
\spline(70,40)(70,50)(80,55)(85,55)(90,50)
\multiput(130,70)(0,-3){6}{\line(0,-1){2}}
\multiput(130,47)(0,-3){13}{\line(0,-1){2}}
\multiput(90,150)(0,-3){6}{\line(0,-1){2}}
\multiput(90,127)(0,-3){13}{\line(0,-1){2}}
\multiput(58,230)(0,-3){6}{\line(0,-1){2}}
\multiput(58,207)(0,-3){13}{\line(0,-1){2}}
\multiput(60,310)(0,-3){6}{\line(0,-1){2}}
\multiput(60,287)(0,-3){13}{\line(0,-1){2}}
\put(30,30){\makebox(0,0){$a$}}
\put(22,23){\makebox(0,0){$x$}}
\put(30,10){\makebox(0,0){$b$}}
\put(60,30){\makebox(0,0){$a$}}
\put(52,23){\makebox(0,0){$y$}}
\put(60,10){\makebox(0,0){$b$}}
\put(30,110){\makebox(0,0){$a$}}
\put(22,103){\makebox(0,0){$x$}}
\put(30,90){\makebox(0,0){$b$}}
\put(60,110){\makebox(0,0){$a$}}
\put(52,103){\makebox(0,0){$y$}}
\put(60,90){\makebox(0,0){$b$}}
\put(30,190){\makebox(0,0){$a$}}
\put(22,183){\makebox(0,0){$x$}}
\put(30,170){\makebox(0,0){$b$}}
\put(30,270){\makebox(0,0){$a$}}
\put(22,263){\makebox(0,0){$c$}}
\put(30,250){\makebox(0,0){$b$}}
\put(45,150){\makebox(0,0){$y$}}
\put(135,53){\makebox(0,0){$a$}}
\put(115,44){\makebox(0,0){$x$}}
\put(145,44){\makebox(0,0){$y$}}
\put(135,27){\makebox(0,0){$b$}}
\put(65,293){\makebox(0,0){$a$}}
\put(45,284){\makebox(0,0){$c$}}
\put(75,284){\makebox(0,0){$c$}}
\put(65,267){\makebox(0,0){$b$}}
\put(95,133){\makebox(0,0){$a$}}
\put(75,124){\makebox(0,0){$x$}}
\put(105,124){\makebox(0,0){$y$}}
\put(95,107){\makebox(0,0){$b$}}
\put(80,25){\makebox(0,0){$b$}}
\put(90,25){\makebox(0,0){$z$}}
\put(95,43){\makebox(0,0){$x$}}
\put(104,44){\makebox(0,0){$w$}}
\put(80,47){\makebox(0,0){$a$}}
\put(85,60){\makebox(0,0){$y$}}
\put(120,135){\makebox(0,0){$a$}}
\put(120,105){\makebox(0,0){$b$}}
\put(130,105){\makebox(0,0){$z$}}
\put(135,123){\makebox(0,0){$x$}}
\put(144,124){\makebox(0,0){$w$}}
\put(140,134){\makebox(0,0){$y$}}
\put(50,210){\makebox(0,0){$a$}}
\put(50,190){\makebox(0,0){$b$}}
\put(70,185){\makebox(0,0){$z$}}
\put(75,203){\makebox(0,0){$x$}}
\put(84,204){\makebox(0,0){$w$}}
\put(80,214){\makebox(0,0){$y$}}
\put(95,283){\makebox(0,0){$x$}}
\put(105,293){\makebox(0,0){$y$}}
\put(104,285){\makebox(0,0){$w$}}
\put(100,303){\makebox(0,0){$z$}}
\put(100,257){\makebox(0,0){$z$}}
\put(1,40){\makebox(0,0){$=$}}
\put(1,120){\makebox(0,0){$=$}}
\put(1,200){\makebox(0,0){$=$}}
\put(0,280){\makebox(0,0){$\displaystyle \sum_c $}}
\end{picture}
\end{center}
\caption{Centrality of $p_{a,b}$ (1)}
\label{nondeg2}
\end{figure}
We next compute the product in the reverse order
as in Figure \ref{nondeg2+}.
\unitlength 0.65mm
\begin{figure}[tb]
\begin{center}
\begin{picture}(160,160)
\thicklines
\put(20,20){\circle{20}}
\put(50,20){\circle{20}}
\put(50,100){\circle{20}}
\put(10,20){\vector(1,0){20}}
\put(40,20){\vector(1,0){20}}
\put(40,100){\vector(1,0){20}}
\put(5,40){\line(1,0){60}}
\put(35,120){\line(1,0){30}}
\put(35,55){\makebox(0,0){$1$}}
\put(50,135){\makebox(0,0){$1$}}
\put(70,40){\vector(1,0){10}}
\put(110,40){\vector(1,0){10}}
\put(140,40){\vector(1,0){10}}
\put(70,120){\vector(1,0){10}}
\put(100,120){\vector(1,0){10}}
\put(0,120){\vector(1,0){10}}
\put(10,130){\vector(1,0){10}}
\put(20,30){\vector(-1,0){0}}
\put(50,30){\vector(-1,0){0}}
\put(20,10){\vector(-1,0){0}}
\put(50,10){\vector(-1,0){0}}
\put(50,110){\vector(-1,0){0}}
\put(50,90){\vector(-1,0){0}}
\put(130,50){\vector(1,0){0}}
\put(126.91,30.49){\vector(1,0){0}}
\put(90,130){\vector(1,0){0}}
\put(86.91,110.49){\vector(1,0){0}}
\put(10,140){\vector(0,-1){40}}
\put(80,30){\vector(1,0){0}}
\put(80,50){\vector(-1,0){0}}
\put(90,30){\vector(1,0){0}}
\put(90,50){\vector(1,0){0}}
\put(75,55){\vector(-1,0){0}}
\put(90,120){\arc{20}{0}{1.157}}
\put(90,120){\arc{20}{1.885}{6.2832}}
\put(130,40){\arc{20}{0}{1.157}}
\put(130,40){\arc{20}{1.885}{6.2832}}
\put(80,40){\arc{20}{0}{4.9}}
\put(80,40){\arc{20}{5.5}{6.2832}}
\put(90,40){\arc{20}{0}{1.9}}
\put(90,40){\arc{20}{2.4}{6.2832}}
\spline(80,30)(67,30)(66,40)(67,55)(80,55)(90,50)
\multiput(130,70)(0,-3){6}{\line(0,-1){2}}
\multiput(130,47)(0,-3){13}{\line(0,-1){2}}
\multiput(90,150)(0,-3){6}{\line(0,-1){2}}
\multiput(90,127)(0,-3){13}{\line(0,-1){2}}
\put(30,30){\makebox(0,0){$a$}}
\put(22,23){\makebox(0,0){$x$}}
\put(30,10){\makebox(0,0){$b$}}
\put(60,30){\makebox(0,0){$a$}}
\put(52,23){\makebox(0,0){$y$}}
\put(60,10){\makebox(0,0){$b$}}
\put(60,110){\makebox(0,0){$a$}}
\put(52,103){\makebox(0,0){$c$}}
\put(60,90){\makebox(0,0){$b$}}
\put(135,53){\makebox(0,0){$a$}}
\put(115,44){\makebox(0,0){$x$}}
\put(145,44){\makebox(0,0){$y$}}
\put(135,27){\makebox(0,0){$b$}}
\put(95,133){\makebox(0,0){$a$}}
\put(75,124){\makebox(0,0){$c$}}
\put(105,124){\makebox(0,0){$c$}}
\put(95,107){\makebox(0,0){$b$}}
\put(82,26){\makebox(0,0){$z$}}
\put(92,26){\makebox(0,0){$b$}}
\put(75,43){\makebox(0,0){$y$}}
\put(76,34){\makebox(0,0){$w$}}
\put(84,43){\makebox(0,0){$a$}}
\put(75,60){\makebox(0,0){$x$}}
\put(5,123){\makebox(0,0){$x$}}
\put(15,134){\makebox(0,0){$y$}}
\put(14,125){\makebox(0,0){$w$}}
\put(10,143){\makebox(0,0){$z$}}
\put(10,97){\makebox(0,0){$z$}}
\put(1,40){\makebox(0,0){$=$}}
\put(24,120){\makebox(0,0){$\displaystyle \sum_c $}}
\end{picture}
\end{center}
\caption{Centrality of $p_{a,b}$ (2)}
\label{nondeg2+}
\end{figure}
Then the coefficients in Figures \ref{nondeg2}
and \ref{nondeg2+} turn out to be
the same, so we have the desired centrality.

The proof of $\sum_{a,b\in{\cal M}}p_{a,b}=1$.
is given graphically in Figure \ref{nondeg3}.
\unitlength 0.65mm
\begin{figure}[tb]
\begin{center}
\begin{picture}(150,130)
\thicklines
\put(50,130){\circle{20}}
\put(80,130){\circle{20}}
\put(65,90){\circle{20}}
\put(50,50){\circle{20}}
\put(35,110){\line(1,0){60}}
\put(40,20){\vector(1,0){20}}
\put(70,50){\vector(1,0){10}}
\put(100,50){\vector(1,0){10}}
\put(100,110){\vector(1,0){10}}
\put(130,110){\vector(1,0){10}}
\put(55,90){\vector(1,0){20}}
\put(50,140){\vector(-1,0){0}}
\put(80,140){\vector(-1,0){0}}
\put(65,80){\vector(-1,0){0}}
\put(65,90){\vector(1,0){0}}
\put(65,100){\vector(-1,0){0}}
\put(50,60){\vector(-1,0){0}}
\put(120,120){\vector(1,0){0}}
\put(116.91,100.49){\vector(1,0){0}}
\put(90,60){\vector(1,0){0}}
\put(88.837,43.343){\vector(1,0){0}}
\put(120,110){\arc{20}{0}{1.157}}
\put(120,110){\arc{20}{1.885}{6.2832}}
\put(90,50){\arc{14}{0}{1.157}}
\put(90,50){\arc{14}{1.885}{6.2832}}
\put(90,50){\arc{20}{3.1416}{6.2832}}
\put(60,140){\makebox(0,0){$a$}}
\put(90,140){\makebox(0,0){$b$}}
\put(75,100){\makebox(0,0){$a$}}
\put(67,93){\makebox(0,0){$c$}}
\put(75,80){\makebox(0,0){$b$}}
\put(125,123){\makebox(0,0){$a$}}
\put(105,114){\makebox(0,0){$c$}}
\put(135,114){\makebox(0,0){$c$}}
\put(125,97){\makebox(0,0){$b$}}
\put(60,60){\makebox(0,0){$b$}}
\put(80,40){\makebox(0,0){$b$}}
\put(105,55){\makebox(0,0){$c$}}
\put(50,25){\makebox(0,0){$c$}}
\multiput(120,140)(0,-3){6}{\line(0,-1){2}}
\multiput(120,117)(0,-3){13}{\line(0,-1){2}}
\multiput(90,80)(0,-3){6}{\line(0,-1){2}}
\multiput(90,54)(0,-3){12}{\line(0,-1){2}}
\put(17,110){\makebox(0,0){$\displaystyle \sum_{a,b} p_{a,b}=
\sum_{a,b,c}$}}
\put(28,50){\makebox(0,0){$=\displaystyle\sum_{b,c}$}}
\put(28,20){\makebox(0,0){$=\displaystyle\sum_c$}}
\put(70,20){\makebox(0,0){$=1$}}
\end{picture}
\end{center}
\caption{$\displaystyle\sum p_{a,b}=1$}
\label{nondeg3}
\end{figure}

We finally have to show that each $p_{a,b}$ is a minimal central projection.
By Ocneanu's theorem mentioned above just before the proof, it is enough
to show that the corresponding $M_\infty$-$M_\infty$ bimodules are all
irreducible.
By the proof of Ocneanu's theorem in \cite{O6} (see \cite[Theorem 4.3]{EK3} or
\cite[Theorem 12.28]{EK4}), we know that the $M_\infty$-$M_\infty$ bimodule
corresponding to $p_{a,b}$ decomposes as
$\bigoplus_{c\in{\cal M}} N_{ab}^c B_c$ as an
$M\vee (M'\cap M_\infty)$-$M_\infty$ bimodule after restricting the
left action to $M\vee (M'\cap M_\infty)$, where $B_c$ denotes the
$M\vee (M'\cap M_\infty)$-$M_\infty$ bimodule labeled with $c$.
This shows that the
Jones index of the $M_\infty$-$M_\infty$ bimodule
corresponding to $p_{a,b}$ is $[a][b]$.  Thus {\sl if} all
the bimodules corresponding to $p_{a,b}$ are irreducible, we get
the global index equal to $\sum_{a,b\in{\cal M}} [a][b]$, which is the
correct global index.  If one or more of the bimodules is reducible, we would
get a smaller global index, which is impossible.  Thus we conclude that
all the bimodules are irreducible.
\qed\end{proof}

We now recall that the
principal graph of the asymptotic inclusion $M\vee (M'\cap M_\infty)
\subset M_\infty$
is given by the fusion graph of the original system
${\cal M}$ by Ocneanu's theorem in \cite{O6}
(see \cite[Theorem 4.1]{EK3} or \cite[Theorem 12.25]{EK4}).
That is,
the set of the odd vertices of the principal graph
is labeled with ${\cal M}$, the
set of the even vertices is labeled by pairs $(a,b)$ with 
$a,b\in {\cal M}$, and the number of the edges between the odd vertex
labeled with $c$ and the even vertex labeled with $(a,b)$ is
given by $N_{ab}^c$, the multiplicity of $c$ in the relative
tensor product $a\otimes_M b$.  The connected component of this
graph containing the even vertices labeled with $(*,*)$ is called the
{\sl fusion graph} of the system ${\cal M}$.  (See \cite{O2},
\cite[page 220]{EK3}, \cite[Section 12.6]{EK4}.)

Combining these pieces of information, 
we get the following proposition.

\begin{proposition}
\label{edges}
Let $N\subset M$ be a hyperfinite type II$_1$ subfactor
with finite index and finite depth.  Suppose that the
system of the $M$-$M$ bimodules arising from this subfactor
has a non-degenerate
braiding.  Then the dual principal graph
of the asymptotic inclusion
$M\vee (M'\cap M_\infty)\subset M_\infty$ is the
fusion graph of the original system, the same as the principal graph.
\end{proposition}

\begin{proof}
As above, we know that the even vertices of the dual principal graph
is labeled with pairs $(a,b)$ for $a,b\in{\cal M}$ and the odd
vertices with $c\in{\cal M}$.  It is thus enough to show that the number
of the edges connecting the vertices labeled with $(a,b)$ and $c$ is
indeed $N_{ab}^c$.  This follows from the above proof of Theorem \ref{ocn}.
\qed\end{proof}

\section{Braiding for $SU(n)_k$ and a tube algebra}
\label{sect-braid}

We now work on the WZW-model $SU(n)_k$.
Let $N\subset M$ be the corresponding Wenzl subfactor in \cite{We}
constructed as in \cite[Section 4]{BG}.
Note that the fusion rule algebra
for the WZW-model $SU(n)_k$ has a natural $\Z/n\Z$-grading and
that the fusion rule subalgebra given by the grade 0 elements
corresponds to the fusion rule algebra of the $M$-$M$ bimodules
arising from this subfactor $N\subset M$.
(This correspondence also follows from \cite[Section 4]{BG}.)
We denote the grading of a primary field $a$ in the model $SU(n)_k$ by
$\gr(a)\in \Z/n\Z$.
Then this system is often degenerate in the sense of the previous section.
Our next aim is to study the asymptotic inclusions for these degenerate cases.
Some statements in this Section hold for a general RCFT in the
sense of \cite{MS} rather than for the WZW-models,
so we make a general statement in such a case.

First we have the following general proposition.

\begin{proposition}
Let ${\cal M}$ be a braided system of $M$-$M$ bimodules.
A bimodule $x\in {\cal M}$ is degenerate if and only if the bimodule
$x$ satisfies the equality in Figure \ref{deg=} for all $y\in{\cal M}$.
\thicklines
\unitlength 1.5mm
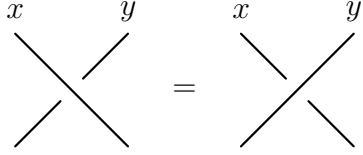
\begin{figure}[tb]
\begin{center}
\begin{picture}(50,25)
\path(10,15)(20,5)
\path(10,5)(14,9)
\path(16,11)(20,15)
\path(30,15)(34,11)
\path(36,9)(40,5)
\path(30,5)(40,15)
\put(10,17){\makebox(0,0){$x$}}
\put(20,17){\makebox(0,0){$y$}}
\put(30,17){\makebox(0,0){$x$}}
\put(40,17){\makebox(0,0){$y$}}
\put(25,10){\makebox(0,0){$=$}}
\end{picture}
\end{center} 
\caption{The degeneracy condition}
\label{deg=}
\end{figure}
\end{proposition}

\begin{proof}
It is trivial that if we have the equality in Figure \ref{deg=},
then we have the degeneracy condition in Figure \ref{degen}.

For the converse direction, we use a graphical argument as in
Figure \ref{deg=deg}, where we have used a handle slide against 
a killing ring.
\thicklines
\unitlength 1.5mm
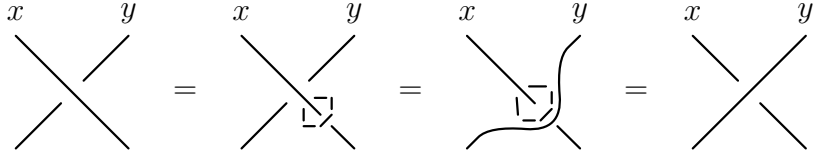
\begin{figure}[tb]
\begin{center}
\begin{picture}(90,25)
\path(10,15)(20,5)
\path(10,5)(14,9)
\path(16,11)(20,15)
\path(30,5)(34,9)
\path(36,11)(40,15)
\path(70,15)(74,11)
\path(76,9)(80,5)
\path(70,5)(80,15)
\path(50,15)(56,9)
\path(58,7)(60,5)
\path(30,15)(37,8)
\path(38,7)(40,5)
\path(37,7)(38,8)
\path(38,8.5)(38,9.5)
\path(37.5,9.5)(36.5,9.5)
\path(35.5,8.5)(35.5,7.5)
\path(35.5,7)(36.5,7)
\path(54.5,8)(54.4,9.5)
\path(55.5,10.5)(57,10.5)
\path(55,7.5)(56,7.5)
\path(57.5,9)(57.5,10)
\path(56.5,7.5)(57.5,8.5)
\spline(50,5)(52,7)(58.5,6.5)(58,13)(60,15)
\put(10,17){\makebox(0,0){$x$}}
\put(20,17){\makebox(0,0){$y$}}
\put(30,17){\makebox(0,0){$x$}}
\put(40,17){\makebox(0,0){$y$}}
\put(50,17){\makebox(0,0){$x$}}
\put(60,17){\makebox(0,0){$y$}}
\put(70,17){\makebox(0,0){$x$}}
\put(80,17){\makebox(0,0){$y$}}
\put(25,10){\makebox(0,0){$=$}}
\put(45,10){\makebox(0,0){$=$}}
\put(65,10){\makebox(0,0){$=$}}
\end{picture}
\end{center} 
\caption{The converse direction}
\label{deg=deg}
\end{figure}
\qed\end{proof}

We record the following straightforward Lemma just to fix the normalization
constants for an RCFT.

\begin{lemma}
The number represented by Figure \ref{Hopf} is $S_{xy}/S_{00}$, where
$S$ denotes the $S$-matrix of the RCFT.
\thicklines
\unitlength 1.2mm
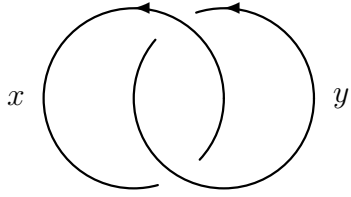
\begin{figure}[tb]
\begin{center}
\begin{picture}(50,30)
\put(30,15){\arc{20}{0}{3.85}}
\put(30,15){\arc{20}{4.4}{6.2832}}
\put(20,15){\arc{20}{0}{0.75}}
\put(20,15){\arc{20}{1.3}{6.2832}}
\put(20,25){\vector(-1,0){0}}
\put(30,25){\vector(-1,0){0}}
\put(7,15){\makebox(0,0){$x$}}
\put(43,15){\makebox(0,0){$y$}}
\end{picture}
\end{center} 
\caption{The Hopf link}
\label{Hopf}
\end{figure}
\end{lemma}

\begin{proof}
This is standard.  See \cite{Wi2} for example.
\qed\end{proof}

Let $\cal M$ be a subsystem of an RCFT.  (A typical case will be
the subsystem of all grade 0 elements in a WZW-model $SU(n)_k$.)
We first have the following lemma.

\begin{lemma}
\label{deg-S}
Let $x$ be an element of the subsystem $\cal \M$.  If we have
$S_{0y}=S_{xy}$ for all $y\in{\cal M}$, then $x$ is degenerate in $\cal M$.
\end{lemma}

\begin{proof}
This follows from a graphical argument as in Figure \ref{deg-S-p}.
\thicklines
\unitlength 1mm
\begin{figure}[tb]
\begin{center}
\begin{picture}(70,70)
\put(25,10){\line(0,1){20}}
\put(45,10){\line(0,1){20}}
\put(60,40){\line(0,1){20}}
\put(20,40){\line(0,1){2}}
\put(20,46){\line(0,1){14}}
\put(20,50){\arc{12}{0}{4.3}}
\put(20,50){\arc{12}{4.95}{6.2832}}
\put(5,50){\makebox(0,0){$\displaystyle \sum_y
\frac{[y]^{1/2}}{[{\cal M}]}$}}
\put(42,50){\makebox(0,0){$\displaystyle =\sum_y
\frac{[y]^{1/2} S_{xy}}{[{\cal M}] S_{x0}}$}}
\put(6,20){\makebox(0,0){$\displaystyle =\sum_y
\frac{[y]^{1/2} S_{0y}}{[{\cal M}] S_{00}}$}}
\put(35,20){\makebox(0,0){$=$}}
\put(25,33){\makebox(0,0){$x$}}
\put(45,33){\makebox(0,0){$x$}}
\put(20,63){\makebox(0,0){$x$}}
\put(60,63){\makebox(0,0){$x$}}
\put(27,55){\makebox(0,0){$y$}}
\end{picture}
\end{center} 
\caption{Degeneracy of $x$}
\label{deg-S-p}
\end{figure}
\qed\end{proof}

\begin{lemma}
\label{deg2}
Suppose that $x$ is degenerate in $\cal M$.  Then for all $y\in{\cal M}$,
we have
\begin{equation}
\label{deg-eq}
\frac{S_{xy}}{S_{00}}=\frac{S_{x0}}{S_{00}}\frac{S_{y0}}{S_{00}}.
\end{equation}
\end{lemma}

\begin{proof}
Suppose that the identity (\ref{deg-eq}) fails for some $y\in{\cal M}$.
Then we have the graphical relation of Figure \ref{deg-neq}.
This, together with identity of Figure \ref{deg-eq-2} given by the
handle slide, gives the identity of Figure \ref{deg-eq-3}, which is
a contradiction.
\thicklines
\unitlength 1mm
\begin{figure}[tb]
\begin{center}
\begin{picture}(105,30)
\put(85,15){\arc{20}{0}{3.85}}
\put(85,15){\arc{20}{4.4}{6.2832}}
\put(75,15){\arc{20}{0}{0.75}}
\put(75,15){\arc{20}{1.3}{6.2832}}
\put(75,25){\vector(-1,0){0}}
\put(85,25){\vector(-1,0){0}}
\put(45,25){\vector(-1,0){0}}
\put(20,25){\vector(-1,0){0}}
\put(64,22){\makebox(0,0){$x$}}
\put(96,22){\makebox(0,0){$y$}}
\put(9,22){\makebox(0,0){$x$}}
\put(56,22){\makebox(0,0){$y$}}
\put(20,15){\circle{20}}
\put(45,15){\circle{20}}
\put(60,15){\makebox(0,0){$\neq$}}
\end{picture}
\end{center} 
\caption{}
\label{deg-neq}
\end{figure}
\thicklines
\unitlength 1mm
\begin{figure}[tb]
\begin{center}
\begin{picture}(85,70)
\put(15,35){\circle{20}}
\put(40,65){\line(0,-1){38}}
\put(40,23){\vector(0,-1){18}}
\put(70,65){\line(0,-1){23}}
\put(70,38){\line(0,-1){26}}
\put(70,8){\vector(0,-1){3}}
\put(70,50){\arc{20}{0}{4.3}}
\put(70,50){\arc{20}{4.95}{6.2832}}
\put(15,45){\vector(-1,0){0}}
\put(70,40){\vector(1,0){0}}
\put(40,35){\arc{20}{0}{0.0374}}
\put(40,35){\arc{20}{0.3366}{0.935}}
\put(40,35){\arc{20}{1.2342}{1.8326}}
\put(40,35){\arc{20}{2.1318}{2.7302}}
\put(40,35){\arc{20}{3.0294}{3.6278}}
\put(40,35){\arc{20}{3.927}{4.5254}}
\put(40,35){\arc{20}{4.8246}{5.423}}
\put(40,35){\arc{20}{5.7222}{6.2832}}
\put(70,20){\arc{20}{0}{0.0374}}
\put(70,20){\arc{20}{0.3366}{0.935}}
\put(70,20){\arc{20}{1.2342}{1.8326}}
\put(70,20){\arc{20}{2.1318}{2.7302}}
\put(70,20){\arc{20}{3.0294}{3.6278}}
\put(70,20){\arc{20}{3.927}{4.5254}}
\put(70,20){\arc{20}{4.8246}{5.423}}
\put(70,20){\arc{20}{5.7222}{6.2832}}
\put(40,68){\makebox(0,0){$x$}}
\put(15,48){\makebox(0,0){$y$}}
\put(70,68){\makebox(0,0){$x$}}
\put(83,50){\makebox(0,0){$y$}}
\put(55,35){\makebox(0,0){$=$}}
\end{picture}
\end{center} 
\caption{}
\label{deg-eq-2}
\end{figure}
\thicklines
\unitlength 1mm
\begin{figure}[tb]
\begin{center}
\begin{picture}(45,55)
\put(20,45){\line(0,-1){28}}
\put(20,13){\vector(0,-1){8}}
\put(20,25){\arc{20}{0}{0.0374}}
\put(20,25){\arc{20}{0.3366}{0.935}}
\put(20,25){\arc{20}{1.2342}{1.8326}}
\put(20,25){\arc{20}{2.1318}{2.7302}}
\put(20,25){\arc{20}{3.0294}{3.6278}}
\put(20,25){\arc{20}{3.927}{4.5254}}
\put(20,25){\arc{20}{4.8246}{5.423}}
\put(20,25){\arc{20}{5.7222}{6.2832}}
\put(20,48){\makebox(0,0){$x$}}
\put(38,25){\makebox(0,0){$=0$}}
\end{picture}
\end{center} 
\caption{}
\label{deg-eq-3}
\end{figure}
\qed\end{proof}

In the rest of this Section, we work on the WZW-model $SU(n)_k$
with $n\mid k$, because it will turn out that this case is a typical
degenerate case related to the orbifold construction.
Let ${\cal M}$ be the subsystem of the WZW-model $SU(n)_k$ consisting of
the elements with grading 0.
(Note that if $n$ and $k$ are relatively prime, then $\cal M$
is non-degenerate by \cite[Section 2]{KT}.)
In this case, the subsystem
$\{x\in{\cal M}\mid S_{0x}=S_{00}\}$ of $\cal M$
is isomorphic to $\Z/n\Z$.  (They are called simple
currents.  See \cite[Section 8.8]{EK4},
\cite[pages 327, 365]{F} for example.)  We
choose and fix an element $\si$ in this subsystem of $\cal M$ so that
this subsystem is given as $\{0,\si,\si^2,\dots,\si^{n-1}\}$.

\begin{lemma}
\label{S-mat1}
For $\sigma$ as above and an arbitrary $y\in{\cal M}$, we have
$S_{0y}=S_{\sigma y}$.
\end{lemma}

\begin{proof}
This follows from  a standard property of the $S$-matrix.  See
\cite[(5.5.25)]{F}, for example.
\qed\end{proof}

This Lemma, together with Lemma \ref{deg-S}, shows that 
any element in
$$\{0,\si,\si^2,\dots,\si^{n-1}\}$$ is degenerate in $\cal M$.
We next show the converse as follows.

\begin{proposition}
\label{WZW-deg}
If $x\in{\cal M}$ is degenerate in $\cal M$, then 
$$x\in \{0,\si,\si^2,\dots,\si^{n-1}\}.$$
\end{proposition}

\begin{proof}
For $y\in \cal M$,
let $\Gamma_y$ be the matrix for the multiplication by $y$ on the entire
fusion rule algebra of the model $SU(n)_k$.  That is,
each entry $(\Gamma_y)_{ab}$ is given by $N_{ay}^b$ for any primary field
$a,b$ in the model $SU(n)_k$.
We also define a vector $v$ by $v=(S_{yx})_y$ for any primary field
$y$ in the model $SU(n)_k$.
According to the grading of $y$, we split the vector $v$ into $n$ pieces
and write $v=(v_0, v_1, \dots, v_{n-1})$, where $v_j$ denotes the vector
component corresponding to $y$ with $\gr(y)=j$.  By the
Verlinde identity \cite{V}, \cite[Section 8.6]{EK4}, we get
$$\Gamma_z v_j =\frac{S_{zx}}{S_{0x}}v_{j+1},$$
where $z\in{\cal M}$ and $j\in \Z/n\Z$.

Lemma \ref{deg2} implies that
we have $S_{zx}\neq 0$ for any $z\in{\cal M}$ with grading 1.
Then we get
\begin{eqnarray*}
\frac{S_{zx}}{S_{0x}}\| v_{j+1}\|_2^2
&=& (\Gamma_z v_j, v_{j+1})\\
&=& (v_j, \Gamma_{\bar z} v_{j+1})\\
&=& \frac{\overline{S_{\bar z x}}}{S_{0x}}\| v_j\|_2^2\\
&=& \frac{S_{z x}}{S_{0x}}\| v_j\|_2^2,
\end{eqnarray*}
which implies $\|v_j\|_2=\| v_{j+1}\|_2$.  Since this is true for all
$j\in \Z/n\Z$ and the matrix $S$ is unitary, we get $\|v_j\|_2=1/\sqrt n$
for all $j\in \Z/n\Z$.  Let $w_0$ be the vector defined by
$(w_0)_y=S_{yx}$
for $y\in{\cal M}$.
Lemma \ref{deg2} implies that $w_0=v_0$ and thus $S_{0x}=S_{00}$.
This means that the Perron--Frobenius weight of the element $x$
is 1 and thus $x$ is in $\{0,\si,\si^2,\dots,\si^{n-1}\}$, which is the
conclusion of the Proposition.
\qed\end{proof}

We now extend the definition of the Ocneanu projection $p_{a,b}$ in
Figure \ref{O-proj}.  Suppose that $a,b$ are primary fields in the model
$SU(n)_k$ with $\gr(a)+\gr(b)=0\in\Z/n\Z$.  Then the graphical formula
in Figure \ref{O-proj} still defines an element in the tube algebra
$\Tube{\cal M}$ for the subsystem $\cal M$ of the elements with
0 grading, because
we have $\gr(c)=0$ for any $c$ appearing in Figure \ref{O-proj}.
Note that $p_{a,b}$ may not be a projection any more.  We call this
element $p_{a,b}$ an {\sl Ocneanu element}.

\begin{lemma}
\label{proj-id5}
For primary fields $a,b$ as above, the element
$p_{a,b}$ is central in the tube algebra $\Tube{\cal M}$.
\end{lemma}

\begin{proof}
The same argument as in Figure \ref{nondeg2} works.
\qed\end{proof}

\begin{lemma}
\label{proj-id}
For primary fields $a,b$ as above, we have
$p_{a,b}=p_{\sigma a, \sigma^{n-1}b}$ in the tube algebra $\Tube{\cal M}$.
\end{lemma}

\begin{proof}
First consider $p_{\sigma, \sigma^{n-1}}$.  
If $c\neq0$ in Figure \ref{O-proj}, then the term corresponding to $c$
is 0.  So we have a single term for this Ocneanu element.
The degeneracy of $\sigma$, proved in Lemmas \ref{S-mat1} and \ref{deg-S},
easily implies $p_{\sigma, \sigma^{n-1}}=p_{0,0}$.

Next note that we have $S(x_1*x_2)S(x_3)=S(x_1)S(x_2*x_3)$ for
$x_1, x_2, x_3 \in H_{S^1\times S^1}$ as in  \cite[Theorem 5.1]{EK3},
\cite[Theorem 12.29]{EK4},
where $S$ means the action of the $S$-matrix in $PSL(2,\Z)$.
(This is a direct analogue of the Verlinde formula \cite{V}.
Actually, the formula in Theorem 5.1 in \cite{EK3} is slightly
incorrect because normalizing coefficients are missing there.
Theorem 12.29 in \cite{EK4} is correct.)
By Lemma \ref{proj-id5}, we can apply this identity to the Ocneanu
elements.  Then we have
$$p_{a,b}=p_{a,b}*p_{0,0}=p_{a,b}*p_{\sigma a, \sigma^{n-1}b}
=p_{\sigma a, \sigma^{n-1}b}.$$
\qed\end{proof}

\begin{lemma}
\label{proj-id3}
Let $a,b,a',b'$ be primary fields in the model $SU(n)_k$ with
$\gr(a)+\gr(b)=\gr(a')+\gr(b')=0\in\Z/n\Z$.  We suppose that
$(a',b')\neq(\sigma^j a,\sigma^{n-j} b)$ for all $j\in\Z/n\Z$.
Then we have $p_{a,b} p_{a', b'}=0$.
\end{lemma}

\begin{proof}
We compute $p_{a,b} p_{a', b'}$ as in Figure \ref{nondeg1}.
The computation is the same up to the third line of Figure \ref{nondeg1}.
Then in the third line, the picture represents the value 0 for any
choice of $d$.  Thus we have $p_{a,b} p_{a', b'}=0$.
\qed\end{proof}

Note that we have a unique primary field
$f$ with $\sigma f=f$, because $n\mid k$.

\begin{lemma}
\label{proj-id6}
If primary fields $a,b$ as above satisfy $(a,b)\neq(f,f)$, then
the element $p_{a,b}$ is a projection in the tube algebra $\Tube{\cal M}$.

Let $P=\{ p_{a,b}\mid \gr(a)+\gr(b)=0\in\Z/n\Z, (a,b)\neq(f,f)\}$.
Then we have $\sum_{p\in P}p+p_{f,f}/n=1$, which implies that
$p_{f,f}/n$ is a central projection.
\end{lemma}

\begin{proof}
Suppose $(a,b)\neq(f,f)$.
We compute $p_{a,b}^2$ as in Figure \ref{nondeg1}.
Then in the third line, we have only the terms with $d$ in
$\{0,\si,\si^2,\dots,\si^{n-1}\}$.  Since $(a,b)\neq(f,f)$,
none of these $d$, except for $d=0$, give a non-zero value.
For $d=0$, we have the original $p_{a,b}$.  This shows that
$p_{a,b}$ is a projection.

Set $P_0=\{ p_{a,b}\mid \gr(a)+\gr(b)=0\in\Z/n\Z\}$.
We compute $\sum_{p \in P_0} p$ as in Figure \ref{proj-id8}.
The second equality follows since the entire system of
the primary fields in the model $SU(n)_k$
is non-degenerate, which follows from unitarity of the $S$-matrix.
(The coefficient $n$
comes from the ratio of the global indices of $\cal M$
and the entire system.)
This implies $\sum_{p\in P}p+p_{f,f}/n=1$.
\unitlength 0.65mm
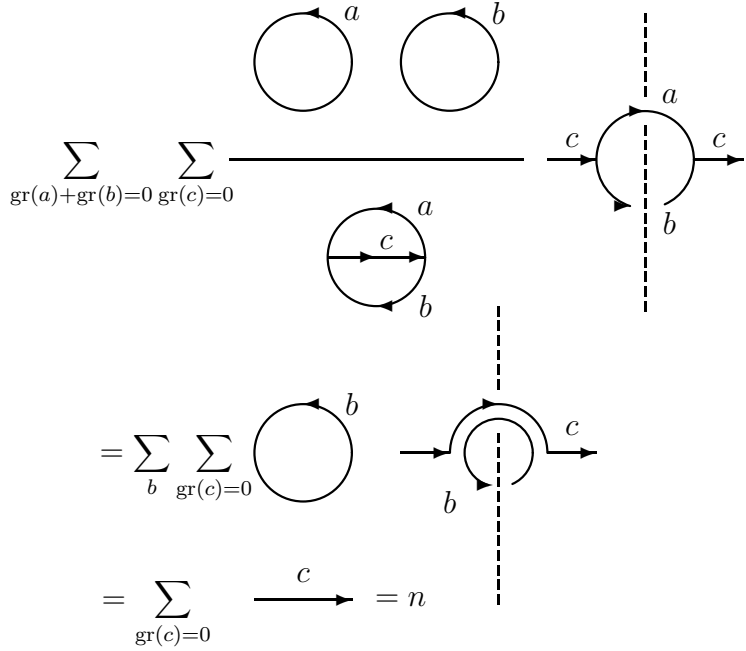
\begin{figure}[tb]
\begin{center}
\begin{picture}(140,130)
\thicklines
\put(40,130){\circle{20}}
\put(70,130){\circle{20}}
\put(55,90){\circle{20}}
\put(40,50){\circle{20}}
\put(25,110){\line(1,0){60}}
\put(30,20){\vector(1,0){20}}
\put(60,50){\vector(1,0){10}}
\put(90,50){\vector(1,0){10}}
\put(90,110){\vector(1,0){10}}
\put(120,110){\vector(1,0){10}}
\put(45,90){\vector(1,0){20}}
\put(40,140){\vector(-1,0){0}}
\put(70,140){\vector(-1,0){0}}
\put(55,80){\vector(-1,0){0}}
\put(55,90){\vector(1,0){0}}
\put(55,100){\vector(-1,0){0}}
\put(40,60){\vector(-1,0){0}}
\put(110,120){\vector(1,0){0}}
\put(106.91,100.49){\vector(1,0){0}}
\put(80,60){\vector(1,0){0}}
\put(78.837,43.343){\vector(1,0){0}}
\put(110,110){\arc{20}{0}{1.157}}
\put(110,110){\arc{20}{1.885}{6.2832}}
\put(80,50){\arc{14}{0}{1.157}}
\put(80,50){\arc{14}{1.885}{6.2832}}
\put(80,50){\arc{20}{3.1416}{6.2832}}
\put(50,140){\makebox(0,0){$a$}}
\put(80,140){\makebox(0,0){$b$}}
\put(65,100){\makebox(0,0){$a$}}
\put(57,93){\makebox(0,0){$c$}}
\put(65,80){\makebox(0,0){$b$}}
\put(115,123){\makebox(0,0){$a$}}
\put(95,114){\makebox(0,0){$c$}}
\put(125,114){\makebox(0,0){$c$}}
\put(115,97){\makebox(0,0){$b$}}
\put(50,60){\makebox(0,0){$b$}}
\put(70,40){\makebox(0,0){$b$}}
\put(95,55){\makebox(0,0){$c$}}
\put(40,25){\makebox(0,0){$c$}}
\multiput(110,140)(0,-3){6}{\line(0,-1){2}}
\multiput(110,117)(0,-3){13}{\line(0,-1){2}}
\multiput(80,80)(0,-3){6}{\line(0,-1){2}}
\multiput(80,54)(0,-3){12}{\line(0,-1){2}}
\put(3,107){\makebox(0,0){$\displaystyle
\sum_{\gr(a)+\gr(b)=0} \sum_{\gr(c)=0}$}}
\put(14,47){\makebox(0,0){$=\displaystyle
\sum_b\sum_{\gr(c)=0}$}}
\put(10,17){\makebox(0,0){$=\displaystyle\sum_{\gr(c)=0}$}}
\put(60,20){\makebox(0,0){$=n$}}
\end{picture}
\end{center}
\caption{$\sum_{p \in P_0} p=n$}
\label{proj-id8}
\end{figure}
The last assertion on $p_{f,f}/n$ now follows from Lemmas
\ref{proj-id3}, \ref{proj-id5}, \ref{proj-id6}.
\qed\end{proof}

\begin{lemma}
\label{proj-id9}
If primary fields $a,b$ as above satisfy $(a,b)\neq(f,f)$, then
the projection $p_{a,b}$ in the tube algebra $\Tube{\cal M}$ is
minimal.
\end{lemma}

\begin{proof}
We define $p_{a,b}^{(c,d)}$ as in Figure \ref{proj-id9-}.
\thicklines
\unitlength 0.7mm
\begin{figure}[tb]
\begin{center}
\begin{picture}(60,80)
\thicklines
\put(10,40){\vector(1,0){10}}
\put(40,40){\vector(1,0){10}}
\put(30,50){\vector(1,0){0}}
\put(26.91,30.49){\vector(1,0){0}}
\put(30,40){\arc{20}{0}{1.157}}
\put(30,40){\arc{20}{1.885}{6.2832}}
\put(35,53){\makebox(0,0){$a$}}
\put(15,43){\makebox(0,0){$c$}}
\put(45,43){\makebox(0,0){$d$}}
\put(35,27){\makebox(0,0){$b$}}
\multiput(30,70)(0,-3){6}{\line(0,-1){2}}
\multiput(30,47)(0,-3){13}{\line(0,-1){2}}
\end{picture}
\end{center} 
\caption{Element $p_{a,b}^{(c,d)}$}
\label{proj-id9-}
\end{figure}
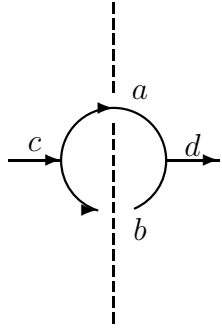

We compute $p_{a,b}^{(c,d)}p_{a,b}^{(d,e)}$
graphically as in Figure \ref{proj-id9+}, where we have
used $(a,b)\neq(f,f)$ in the second line.
\thicklines
\unitlength 0.7mm
\begin{figure}[tb]
\begin{center}
\begin{picture}(140,240)
\put(20,20){\circle{20}}
\put(50,20){\circle{20}}
\put(35,60){\circle{20}}
\put(20,100){\circle{20}}
\put(50,100){\circle{20}}
\put(35,140){\circle{20}}
\put(5,40){\line(1,0){60}}
\put(5,120){\line(1,0){60}}
\put(25,60){\vector(1,0){20}}
\put(25,140){\vector(1,0){20}}
\put(20,30){\vector(-1,0){0}}
\put(50,30){\vector(-1,0){0}}
\put(35,50){\vector(-1,0){0}}
\put(35,70){\vector(-1,0){0}}
\put(20,110){\vector(-1,0){0}}
\put(50,110){\vector(-1,0){0}}
\put(35,130){\vector(-1,0){0}}
\put(35,150){\vector(-1,0){0}}
\put(90,50){\vector(1,0){0}}
\put(86.91,30.49){\vector(1,0){0}}
\put(30,210){\vector(1,0){0}}
\put(26.91,190.49){\vector(1,0){0}}
\put(60,210){\vector(1,0){0}}
\put(56.91,190.49){\vector(1,0){0}}
\put(70,40){\vector(1,0){10}}
\put(100,40){\vector(1,0){10}}
\put(70,120){\vector(1,0){10}}
\put(120,120){\vector(1,0){10}}
\put(100,130){\vector(1,0){0}}
\put(100,110){\vector(1,0){0}}
\put(100,120){\vector(0,-1){0}}
\put(10,200){\vector(1,0){10}}
\put(40,200){\vector(1,0){10}}
\put(70,200){\vector(1,0){10}}
\put(90,130){\line(1,0){20}}
\put(93,110){\line(1,0){17}}
\put(100,110){\line(0,1){5}}
\put(100,120){\line(0,1){10}}
\put(90,40){\arc{20}{0}{1.157}}
\put(90,40){\arc{20}{1.885}{6.2832}}
\put(30,200){\arc{20}{0}{1.157}}
\put(30,200){\arc{20}{1.885}{6.2832}}
\put(60,200){\arc{20}{0}{1.157}}
\put(60,200){\arc{20}{1.885}{6.2832}}
\put(110,120){\arc{20}{0}{1.5708}}
\put(110,120){\arc{20}{4.7124}{6.2832}}
\put(90,120){\arc{20}{1.885}{4.7124}}
\multiput(90,70)(0,-3){6}{\line(0,-1){2}}
\multiput(90,47)(0,-3){13}{\line(0,-1){2}}
\multiput(90,150)(0,-3){6}{\line(0,-1){2}}
\multiput(90,127)(0,-3){13}{\line(0,-1){2}}
\multiput(30,230)(0,-3){6}{\line(0,-1){2}}
\multiput(30,207)(0,-3){13}{\line(0,-1){2}}
\multiput(60,230)(0,-3){6}{\line(0,-1){2}}
\multiput(60,207)(0,-3){13}{\line(0,-1){2}}
\put(99,117){\line(1,0){2}}
\put(102,118){\line(1,1){2}}
\put(97,118){\line(-1,1){2}}
\put(95,121){\line(0,1){2}}
\put(95,125){\line(1,1){2}}
\put(101,127){\line(1,0){2}}
\put(30,10){\makebox(0,0){$a$}}
\put(60,10){\makebox(0,0){$b$}}
\put(45,70){\makebox(0,0){$a$}}
\put(37,63){\makebox(0,0){$d$}}
\put(45,50){\makebox(0,0){$b$}}
\put(30,90){\makebox(0,0){$a$}}
\put(60,90){\makebox(0,0){$b$}}
\put(45,150){\makebox(0,0){$a$}}
\put(37,143){\makebox(0,0){$d$}}
\put(45,130){\makebox(0,0){$b$}}
\put(95,53){\makebox(0,0){$a$}}
\put(75,44){\makebox(0,0){$c$}}
\put(105,44){\makebox(0,0){$e$}}
\put(95,27){\makebox(0,0){$b$}}
\put(35,213){\makebox(0,0){$a$}}
\put(15,204){\makebox(0,0){$c$}}
\put(45,204){\makebox(0,0){$d$}}
\put(35,187){\makebox(0,0){$b$}}
\put(65,213){\makebox(0,0){$a$}}
\put(75,204){\makebox(0,0){$e$}}
\put(65,187){\makebox(0,0){$b$}}
\put(75,124){\makebox(0,0){$c$}}
\put(125,124){\makebox(0,0){$e$}}
\put(82,133){\makebox(0,0){$a$}}
\put(82,107){\makebox(0,0){$b$}}
\put(118,133){\makebox(0,0){$a$}}
\put(118,107){\makebox(0,0){$b$}}
\put(104,123){\makebox(0,0){$0$}}
\put(0,40){\makebox(0,0){$=$}}
\put(0,120){\makebox(0,0){$=$}}
\end{picture}
\end{center} 
\caption{Product $p_{a,b}^{(c,d)}p_{a,b}^{(d,e)}$}
\label{proj-id9+}
\end{figure}
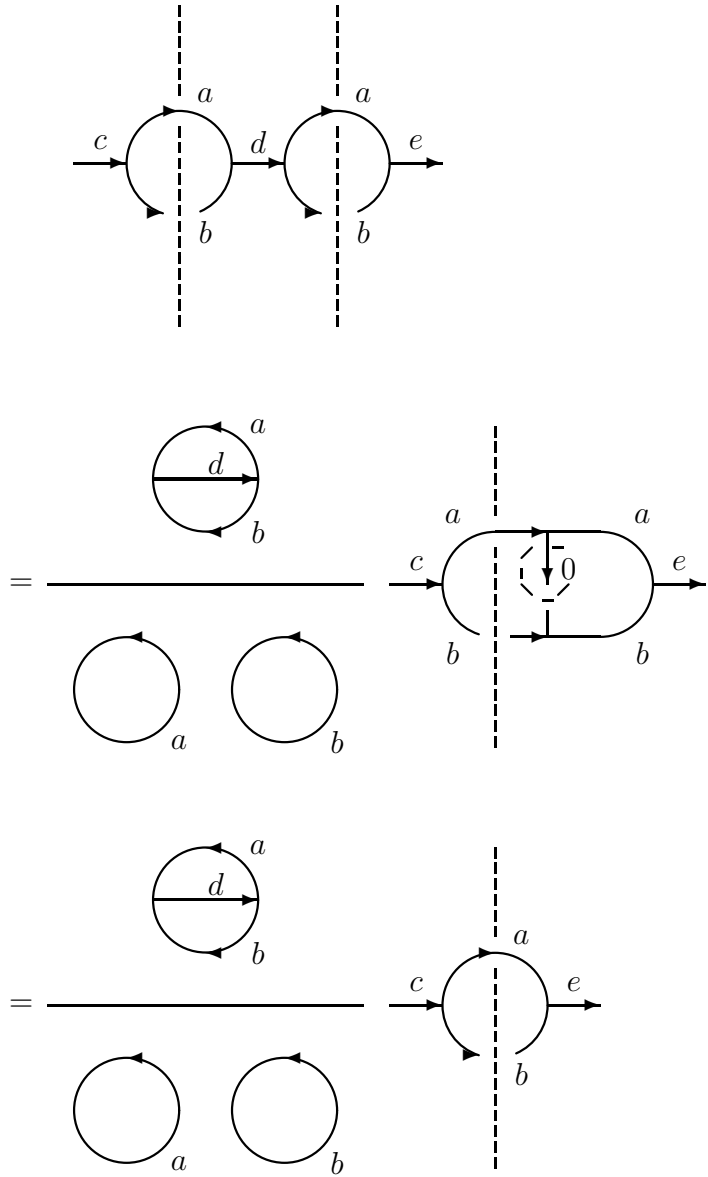
Let $M(a,b)$ be the set of the 
primary fields $c$ satisfying
$p_{a,b}^{(c,c)}\neq 0$.
If $c,d\in M(a,b)$, then 
the computation in Figure \ref{proj-id9+} shows that
$p_{a,b}^{(c,d)}\neq0$.
If $c,d,e\in M(a,b)$, then
the computation in Figure \ref{proj-id9+} also shows that
$p_{a,b}^{(c,d)}p_{a,b}^{(d,e)}\neq0$.
We thus have a system of matrix units
$\{\lambda_{c,d}p_{a,b}^{(c,d)}\}_{c,d}$ for
$p_{a,b}(\Tube {\cal M})p_{a,b}$, where $\lambda_{c,d}$
are some positive numbers.  This shows that the algebra
$p_{a,b}(\Tube {\cal M})p_{a,b}$ is a full matrix algebra
and thus the central projection $p_{a,b}$ is minimal in
the center of the tube algebra $\Tube{\cal M}$.
\qed\end{proof}

Again by Ocneanu's theorem in \cite{O6} (see \cite[Theorem 4.3]{EK3} or
\cite[Theorem 12.28]{EK4}), we get irreducible $M_\infty$-$M_\infty$
bimodules corresponding to $p_{a,b}$ with $(a,b)\neq(f,f)$.
Note that we have corresponding bimodules even when
$\gr(a),\gr(b)\neq0$.  The principal graphs of the asymptotic
inclusions are determined
only by the primary fields with grading 0, but the dual principal
graphs have vertices related to the primary fields with other
grading.  The primary fields with non-zero grading are the
{\sl ghosts} of the system $\cal M$ in the sense of Ocneanu
\cite{O7}.

We work on the irreducible
decompositions of these $M_\infty$-$M_\infty$
bimodules after restricting the left action 
to $M\vee (M'\cap M_\infty)$ as follows.  This gives partial
information about the dual principal graph of the asymptotic
inclusion.

\begin{lemma}
\label{proj-id10}
Let $a,b$ be primary fields as above satisfying $(a,b)\neq(f,f)$
and $X_{a,b}$ the irreducible
$M_\infty$-$M_\infty$ bimodule corresponding to
the minimal central projection $p_{a,b}$ in $\Tube{\cal M}$.
If we restrict the left action  to $M\vee (M'\cap M_\infty)$,
we get a decomposition $X_{a,b}=\bigoplus_{c\in {\cal M}} N_{ab}^c X_c$
as a $M\vee (M'\cap M_\infty)$-$M_\infty$ bimodule, where
$X_c$ is the $M\vee (M'\cap M_\infty)$-$M_\infty$ bimodule
corresponding to $c\in{\cal M}$.
\end{lemma}

\begin{proof}
An argument similar to the one in
the proof of Theorem \ref{ocn} works.
\qed\end{proof}

We next work on the case $(a,b)=(f,f)$.  We still have an
$M_\infty$-$M_\infty$ bimodule in such a case, though this
bimodule might not be irreducible, and we get the following
lemma in the same way.

\begin{lemma}
\label{proj-id11}
Let $X_{f,f}$ be the $M_\infty$-$M_\infty$ bimodule corresponding to
the central projection $p_{f,f}/n$ in $\Tube{\cal M}$.
If we restrict the left action  to $M\vee (M'\cap M_\infty)$,
we get a decomposition $X_{f,f}=\bigoplus_{c\in {\cal M}} N_{ff}^c X_c$
as a $M\vee (M'\cap M_\infty)$-$M_\infty$ bimodule, where
$X_c$ is the $M\vee (M'\cap M_\infty)$-$M_\infty$ bimodule
corresponding to $c\in{\cal M}$.
\end{lemma}

\begin{proof}
An argument similar to the one in
the proof of Theorem \ref{ocn} again works.
\qed\end{proof}

We would like to get a full description of the dual principal
graph, but the bimodule $X_{f,f}$ plays a quite subtle role.
So we first make the following assumption and later prove
that this assumption holds in some cases.

\begin{assumption}
\label{orbif}
The $M_\infty$-$M_\infty$ bimodule $X_{f,f}$ decomposes
into $n$ irreducible bimodules and each has the same dimension.
\end{assumption}

In this assumption, we mean the square root of the Jones
index of the corresponding subfactor of a bimodule
by the ``dimension'' of a bimodule.
A. Ocneanu has observed this assumption holds for $SU(2)_{2k}$ and
we will prove that this also holds for $SU(3)_{3k}$ in a general
framework.  We conjecture that this assumption holds for any
$SU(n)_{nk}$, but combinatorial complexity has prevented us
from proving it, so far.

A simple computation easily gives the following lemma.

\begin{lemma}
\label{glob}
Assumption \ref{orbif} gives the correct global index for the
asymptotic inclusion.
\end{lemma}

Consider the dual principal graph of the asymptotic
inclusion.
To each $M_\infty$-$M_\infty$ or
$M\vee (M'\cap M_\infty)$-$M_\infty$ bimodule, we assign
its dimension, as usual.   This gives a Perron--Frobenius weight.
That is, for an $M\vee (M'\cap M_\infty)$-$M_\infty$ bimodule
corresponding to $c\in{\cal M}$, we get $[c]^{1/2}[{\cal M}]^{1/2}$
and for an $M_\infty$-$M_\infty$ bimodule
corresponding to $p_{a,b}$ with arbitrary
$a,b$ in the model $SU(n)_k$ with $(a,b)\neq(f,f)$,
we get $[a]^{1/2} [b]^{1/2}$.  We also note that the
Perron--Frobenius eigenvalue for this weight is $[{\cal M}]^{1/2}$.

\begin{lemma}
\label{PF}
These Perron--Frobenius weights on the
$M\vee (M'\cap M_\infty)$-$M_\infty$ bimodules
are compatible with Assumption \ref{orbif}.
\end{lemma}

\begin{proof}
For $c\in{\cal M}$, we denote by $X_c$ the corresponding
$M\vee (M'\cap M_\infty)$-$M_\infty$ bimodule.
We easily get $[X_c]=[c][{\cal M}]$.
We can also form a fusion  graph using all the
primary fields in the model $SU(n)_k$.  From the Perron--Frobenius
property of this graph, we get
$$n[{\cal M}] [c]^{1/2}=
\sum_{a,b} N_{ab}^c [a]^{1/2}[b]^{1/2},$$
where $a,b$ are arbitrary primary fields in the model $SU(n)_k$.
Let $L$ be a set of representatives of the equivalence
classes on all the pairs of  arbitrary primary fields in the model
$SU(n)_k$ excluding $(f,f)$ for the equivalence relation
$(a,b)\sim (\sigma^j a,\sigma^{n-j} a)$ with $j\in\Z/n\Z$.
Then the right hand side of the equality is equal to
$$n \left(\sum_{(a,b)\in L} N_{ab}^c [a]^{1/2}[b]^{1/2}+
N_{ff}^c \frac{[f]}{n}\right),$$
by Assumption \ref{orbif}.  This gives
$$[{\cal M}]^{1/2} [c]^{1/2}[{\cal M}]^{1/2}=
\sum_{(a,b)\in L} N_{ab}^c [a]^{1/2}[b]^{1/2}+
N_{ff}^c \frac{[f]}{n},$$
which is the conclusion, because we have
Lemmas \ref{proj-id10}, \ref{proj-id11}.
\qed\end{proof}

\section{Dual principal graphs of the asymptotic inclusions
--- $SU(2)_k$ case ---}
\label{sect-su2}

With the preliminaries of the previous section, we compute
the dual principal graphs of the asymptotic inclusions of
the $SU(2)_{2n}$ subfactors, that is, the Jones subfactors
of type $A_{2n+1}$ constructed in \cite{J}, with $n>1$.
These results were first claimed by Ocneanu.  We present a complete
proof here, because we will generalize the results in the next
Section.

First label the primary fields in $SU(2)_{2n}$ with
$0,1,\dots,2n$ as usual.  Recall that the fusion rule is given as
$$N_{jk}^l=\left\{\begin{array}{ll}
1,&\quad{\rm if \ }|j-k|\le l\le j+k, j+k+l\in 2\Z,
j+k+l\le 4n,\\
0,&\quad{\rm otherwise.}\end{array}\right.$$

Note that all the bimodules in ${\cal M}$ are labeled with
even integers and they are all self-contragredient.
The $M\vee (M'\cap M_\infty)$-$M_\infty$ bimodules
arising from the asymptotic inclusion are labeled
with $0,2,\dots,2n$ and the
$M\vee (M'\cap M_\infty)$-$M\vee (M'\cap M_\infty)$
bimodules are labeled with pairs of even integers
$0,2,\dots,2n$.  This implies that all the
$M\vee (M'\cap M_\infty)$-$M\vee (M'\cap M_\infty)$
bimodules arising from the asymptotic inclusion are
also self-contragredient.

On the even vertices of the dual principal graph,
we do not know how the $M_\infty$-$M_\infty$ bimodule
corresponding to $p_{n,n}/2$ in $\Tube {\cal M}$ decomposes
into irreducible ones, but the fusion rule as above shows
that this bimodule contains exactly one copy of $X_0$
when we restrict the left action to $M\vee (M'\cap M_\infty)$
by Lemma \ref{proj-id11}.  Then Lemma \ref{PF}
implies that the $M_\infty$-$M_\infty$ bimodule
corresponding to $p_{n,n}/2$ is not irreducible and
it contains at least one irreducible bimodule
whose dimension is half the dimension of this 
$M_\infty$-$M_\infty$ bimodule corresponding to $p_{n,n}/2$.
Then Lemma \ref{glob} implies that this 
$M_\infty$-$M_\infty$ bimodule corresponding to $p_{n,n}/2$
decomposes into exactly two irreducible bimodules with
equal Jones indices and thus Assumption \ref{orbif}
holds, because we would have a smaller global index
otherwise.  We label these two bimodules with $(n,n)_+$
and $(n,n)_-$.
We will now determine the dual principal
graph of the asymptotic inclusion.
By Lemma \ref{proj-id10}, it is enough to determine how the
even vertices labeled with $(n,n)_+$ and $(n,n)_-$
are connected to the odd vertices.  Since the odd
vertex labeled with 0 is connected to one of these
two even vertices, we may assume that $(n,n)_+$ is
connected to 0.

\begin{lemma}
\label{self-c}
The $M_\infty$-$M_\infty$ bimodules labeled with
$(n,n)_\pm$ are self-contragredient.
\end{lemma}

\begin{proof}
First note that the other
$M_\infty$-$M_\infty$ bimodules
are self-contragredient by Lemma \ref{proj-id}.

We count the number of the paths of length 2
connecting the odd vertex
$0$ to itself on the principal graph of the asymptotic
inclusion, which is the fusion graph of ${\cal M}$,
via the contragredient map, because the fusion graph is
now connected.
By the fusion rule described above,
we can go from $0$ back to $0$ on the
principal graph through $(0,0),(2,2),\dots,(2n,2n)$.
This implies that the number of the paths is $n+1$.

We know that the number of paths of length 2
connecting the odd vertex
$0$ to itself on the dual principal graph of the asymptotic
inclusion via the contragredient map is also equal to $n-1$ by
(bi)unitarity of the connection arising from the
asymptotic inclusion
\cite[page 130]{O1} (or \cite[Section 10.3]{EK4}).

The $M_\infty$-$M_\infty$ bimodules labeled with
$$(0,0)=(2n,2n), (1,1)=(2n-1,2n-1),\dots,(n-1,n-1)=(n+1,n+1)$$
give $n$ paths from $0$ back to $0$ on the dual principal
graph.  (Here the equality as in $(0,0)=(2n,2n)$ means that
the bimodule labeled with $p_{0,0}$ is equal to that
with $p_{2n,2n}$ because of Lemma \ref{proj-id}.)
This means that we still have another path from $0$
back to $0$ on the dual principal graph through the
even vertex labeled with $(n,n)_+$. 

This means that the 
$M_\infty$-$M_\infty$ bimodule labeled with
$(n,n)_+$, hence that with $(n,n)_-$, is contragredient to itself.
\qed\end{proof}

We next count the number of paths connecting the odd vertices
$0$ and $2$ on the principal graph of the asymptotic
inclusion.  (In this kind of counting in the rest of this
paper, by a ``path'' we mean a path of length 2 on the graph via
the contragredient map.)
Again by the fusion rule, we can go from
$0$ to $2$ on the
principal graph through $(2,2),(4,4),\dots,(2n-2,2n-2)$.
This implies that the number of the paths is $n-1$.

Again by unitarity,
the number of the paths connecting the odd vertices
$0$ and $2$ on the dual principal graph is also equal
to $n-1$.
The even vertices labeled with $(1,1),(2,2),\dots,(n-1,n-1)$
are connected both to the odd vertices $0$ and $2$ by
Lemma \ref{proj-id10}.  These already give the correct
number of paths, so this fact means that the even vertex
$(n,n)_+$ is not connected to the odd
vertex $2$.  Then Lemma \ref{proj-id11} implies that
the even vertex $(n,n)_-$ is connected to the odd
vertex $2$.

Similarly, we can count the number of paths from $0$ to $4,6,\dots$
on the principal/dual principal graphs with the fusion rule.
Then unitarity gives the following description of the
dual principal graph.

\begin{theorem}
\label{dual-su2}
Let $N\subset M$ be the subfactor corresponding to
$SU(2)_{2n}$.  Then the even vertex $(n,n)_+$ of the
dual principal graph of the asymptotic inclusion is
connected to the odd vertices
$0,4,\dots$.  The even vertex $(n,n)_-$ of the
dual principal graph is connected to the odd vertices
$2,6,\dots$.
\end{theorem}

As a corollary of the above description,
we get the following, which was announced
by Ocneanu in \cite[page 41]{O7}.  Note that this Corollary
gives the number of the even vertices of the dual principal
graph of the asymptotic inclusions.  These are also the dimensions
of the Hilbert spaces $H_{S^1\times S^1}$ in the corresponding
topological quantum field theories.

\begin{corollary}
\label{ocn2}
Let $N\subset M$ be the subfactor corresponding to
$SU(2)_k$, that is, the Jones subfactor of type
$A_{k+1}$.  Assume $k>2$.  Then the number of the
irreducible $M_\infty$-$M_\infty$ bimodules arising
from the asymptotic inclusion is given as follows.
$$\begin{array}{ll}
\displaystyle\left(\frac{k+1}{2}\right)^2,
&\quad {\rm if \ }k {\rm \ is \ odd,}\\
\displaystyle\frac{k^2}{4}+\frac{k}{2}+2,
&\quad {\rm if \ }k {\rm \ is \ even.}
\end{array}$$
\end{corollary}

We list some examples of the  dual principal
graphs.  The first one is for $SU(2)_4$, which is the
Jones subfactor of type $A_5$.  It is well known
that this subfactor of index 3 is
of the form $R\rtimes S_2\subset R\rtimes S_3$,
where $S_2$ and $S_3$ are the symmetric groups of order
2 and 3 respectively and these groups act freely on
the hyperfinite II$_1$ factor $R$.  (See \cite{O1}.)
Thus the paragroup
of the asymptotic inclusion is given by that of the
subfactor $R^{S_3\times S_3}\subset R^{S_3}$, where
$S_3$ is diagonally embedded into $S_3\times S_3$ and
the group $S_3$ acts freely on $R$, by Ocneanu's theorem.
(See \cite[Lemma 2.15]{K2}, \cite[Appendix]{K3},
\cite[Section 12.8]{EK4}.)

So the (dual) principal graphs of the asymptotic inclusion
can be described with Ocneanu's theorem on subfactors
of the form $R^G\subset R^H$, where $G$ is a finite group
acting freely on a II$_1$ factor $R$ and $H$ is a subgroup
of $G$.  (See \cite{KY} for this type of computation.)
Of course, this method gives the same result as in
Figure \ref{su24}.

\unitlength 0.8mm
\begin{figure}[tb]
\begin{center}
\begin{picture}(160,70)
\thinlines
\multiput(10,20)(20,0){8}{\circle*{2}}
\multiput(30,60)(50,0){3}{\circle*{2}}
\path(30,60)(10,20)
\path(30,60)(30,20)
\path(30,60)(50,20)
\path(80,60)(30,20)
\path(80,60)(70,20)
\path(80,60)(90,20)
\path(80,60)(110,20)
\path(80,60)(130,20)
\path(130,60)(50,20)
\path(130,60)(130,20)
\path(130,60)(150,20)
\put(10,15){\makebox(0,0){$00$}}
\put(10,8){\makebox(0,0){$44$}}
\put(30,15){\makebox(0,0){$11$}}
\put(30,8){\makebox(0,0){$33$}}
\put(50,15){\makebox(0,0){$22_+$}}
\put(110,15){\makebox(0,0){$22_-$}}
\put(70,15){\makebox(0,0){$02$}}
\put(70,8){\makebox(0,0){$42$}}
\put(90,15){\makebox(0,0){$20$}}
\put(90,8){\makebox(0,0){$24$}}
\put(130,15){\makebox(0,0){$13$}}
\put(130,8){\makebox(0,0){$31$}}
\put(150,15){\makebox(0,0){$04$}}
\put(150,8){\makebox(0,0){$40$}}
\put(30,65){\makebox(0,0){$0$}}
\put(80,65){\makebox(0,0){$2$}}
\put(130,65){\makebox(0,0){$4$}}
\end{picture}
\end{center}
\caption{Dual principal graphs for $SU(2)_4$, $A_5$}
\label{su24}
\end{figure}
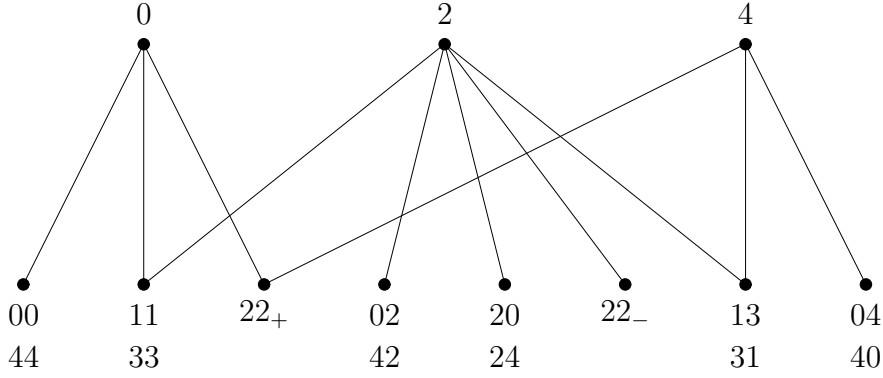

A more complicated example of the dual principal graph
of the asymptotic inclusion is given in Figure \ref{su26}.

\unitlength 0.5mm
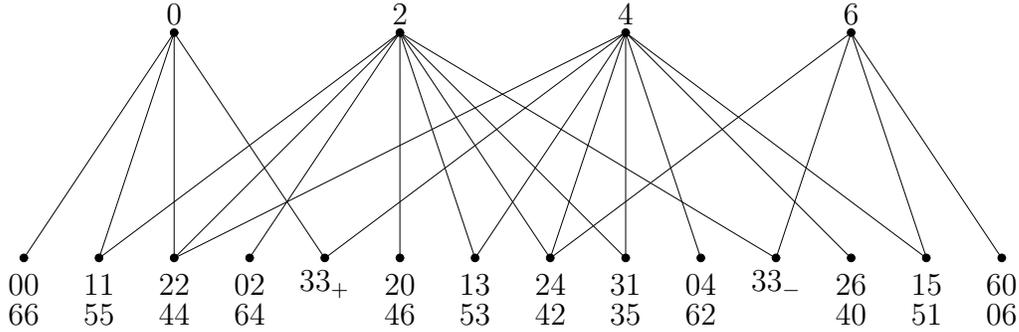
\begin{figure}[tb]
\begin{center}
\begin{picture}(280,90)
\thinlines
\multiput(10,20)(20,0){14}{\circle*{2}}
\multiput(50,80)(60,0){4}{\circle*{2}}
\path(50,80)(10,20)
\path(50,80)(30,20)
\path(50,80)(50,20)
\path(50,80)(90,20)
\path(110,80)(30,20)
\path(110,80)(50,20)
\path(110,80)(70,20)
\path(110,80)(110,20)
\path(110,80)(130,20)
\path(110,80)(150,20)
\path(110,80)(170,20)
\path(110,80)(210,20)
\path(170,80)(50,20)
\path(170,80)(90,20)
\path(170,80)(130,20)
\path(170,80)(150,20)
\path(170,80)(170,20)
\path(170,80)(190,20)
\path(170,80)(230,20)
\path(170,80)(250,20)
\path(230,80)(150,20)
\path(230,80)(210,20)
\path(230,80)(250,20)
\path(230,80)(270,20)
\put(90,13){\makebox(0,0){$33_+$}}
\put(210,13){\makebox(0,0){$33_-$}}
\put(10,13){\makebox(0,0){$00$}}
\put(10,5){\makebox(0,0){$66$}}
\put(30,13){\makebox(0,0){$11$}}
\put(30,5){\makebox(0,0){$55$}}
\put(50,13){\makebox(0,0){$22$}}
\put(50,5){\makebox(0,0){$44$}}
\put(70,13){\makebox(0,0){$02$}}
\put(70,5){\makebox(0,0){$64$}}
\put(110,13){\makebox(0,0){$20$}}
\put(110,5){\makebox(0,0){$46$}}
\put(130,13){\makebox(0,0){$13$}}
\put(130,5){\makebox(0,0){$53$}}
\put(150,13){\makebox(0,0){$24$}}
\put(150,5){\makebox(0,0){$42$}}
\put(170,13){\makebox(0,0){$31$}}
\put(170,5){\makebox(0,0){$35$}}
\put(190,13){\makebox(0,0){$04$}}
\put(190,5){\makebox(0,0){$62$}}
\put(230,13){\makebox(0,0){$26$}}
\put(230,5){\makebox(0,0){$40$}}
\put(250,13){\makebox(0,0){$15$}}
\put(250,5){\makebox(0,0){$51$}}
\put(270,13){\makebox(0,0){$60$}}
\put(270,5){\makebox(0,0){$06$}}
\put(50,85){\makebox(0,0){$0$}}
\put(110,85){\makebox(0,0){$2$}}
\put(170,85){\makebox(0,0){$4$}}
\put(230,85){\makebox(0,0){$6$}}
\end{picture}
\end{center}
\caption{Dual principal graphs for $SU(2)_6$, $A_7$}
\label{su26}
\end{figure}

In the graphs in Figures \ref{su24}, \ref{su26},
the vertices labeled with pairs of odd numbers arise,
while the original $M$-$M$ bimodules are labeled with
only even numbers.  These odd numbers correspond to
the ghosts in the terminology of Ocneanu \cite{O7}.

\section{Dual principal graphs of the asymptotic inclusions
--- $SU(3)_k$ case ---}
\label{sect-su3}

Now we work on the asymptotic inclusions of the
$SU(3)_{3k}$-subfactors and give our main results in this paper.
We have to determine how the central projection
$p_{f,f}/3$ decomposes into minimal central projections
in the tube algebra $\Tube {\cal M}$.

Lemmas \ref{proj-id11} and \ref{PF} imply that $p_{f,f}/3$
contains at least one minimal central projection $p^{(0)}_{f,f}$,
the dimension of the corresponding irreducible
$M_\infty$-$M_\infty$
bimodule of which is one third of that
of the $M_\infty$-$M_\infty$ bimodule corresponding to 
$p_{f,f}/3$.  This argument also shows that
the odd vertex of the dual principal graph labeled with $0$
is connected to the even vertex labeled with $p^{(0)}_{f,f}$ with
exactly one edge and also that the odd vertex $0$
is not connected to the other even vertices arising from
the decomposition of $p_{f,f}/3$.

\begin{lemma}
\label{self-c2}
The irreducible $M_\infty$-$M_\infty$ bimodule corresponding to 
$p^{(0)}_{f,f}$ is contragredient to itself.
\end{lemma}

\begin{proof}
We again count the number of the appropriate paths as in
the proof of Lemma \ref{self-c}.

We first count the number of the paths connecting the odd vertex
$0$ to itself on the principal graph of the asymptotic
inclusion.  Let $l$ be the number of the primary fields in the
WZW-model $SU(3)_{3k}$.  Then it is easy to see that
the number of the
primary fields in ${\cal M}$ is $(l+2)/3$.
By the fusion rule, the number of paths from $0$
back to $0$ on the principal graph is $(l+2)/3$.

It is also easy to see
that the number of paths connecting the odd vertex
$0$ to itself on the dual principal graph of the asymptotic
inclusion without going through the even vertices corresponding
to some minimal central projection appearing in the
decomposition of $p_{f,f}/3$ is $(l-1)/3$.

These mean that we still have one more path from $0$
back to $0$ on the dual principal graph, which must
go through the even vertex corresponding to 
$p^{(0)}_{f,f}$.   That is, the
$M_\infty$-$M_\infty$ bimodule corresponding to 
$p^{(0)}_{f,f}$ is self-contragredient.
\qed\end{proof}

\begin{lemma}
\label{edge1}
If $c\in{\cal M}$ satisfies $N_{ff}^c=1$, then the odd
vertex of the dual principal graph labeled with $c$ is
connected to the even vertex labeled with $p^{(0)}_{f,f}$
with exactly one edge.
\end{lemma}

\begin{proof}
We count the number of appropriate paths again.

The number of paths connecting the odd vertex
$0$ to itself on the principal graph of the asymptotic
inclusion is $\sum_{a\in{\cal M}} N_{a\bar a}^c$, because
the principal graph is the fusion graph which is now connected.

Let $l$ be the number of the edges connecting the odd
vertex of the dual principal graph labeled with $c$ and
the even vertex labeled with $p^{(0)}_{f,f}$.
Lemma \ref{proj-id10} implies that $l$ is 0 or 1, because
$N_{ff}^c=1$.
We next count the number of paths connecting the odd vertex
$0$ to itself on the dual principal graph of the asymptotic
inclusion.  This number is equal to
$(\sum_a N_{a \bar a}^c-N_{ff}^c)/3+l$, where the summation
is over all the primary fields $a$ of the WZW-model $SU(3)_{3k}$.

Since the two numbers are equal, we get
$$\sum_{\gr(a)=1,2} N_{a \bar a}^c=
2\sum_{\gr(a)=1} N_{a \bar a}^c=-1+3l,$$
which implies that $3l-1$ is even.  That is, we get $l=1$.
\qed\end{proof}

As in Section \ref{sect-braid}, we know that
the subsystem
$\{x\in{\cal M}\mid S_{0x}=S_{00}\}$ of $\cal M$
is given as $\{0,\si,\si^2\}$.  The above Lemma gives
the following.

\begin{corollary}
\label{cor1}
Each of the odd
vertices of the dual principal graph labeled with $0,\si,\si^2$
is connected to the even vertex labeled with $p^{(0)}_{f,f}$
with exactly one edge.
\end{corollary}

\begin{proof}
This follows from the above Lemma, because we have
$N_{ab}^c=N_{a\sigma(b)}^{\sigma(c)}$ by \cite{Wal}.
\qed\end{proof}

We now need some lemmas for the fusion rule of the 
WZW-model $SU(3)_{3k}$, which has been obtained by
Goodman--Wenzl in \cite{GW} as a quantum version of the
classical Littlewood--Richardson rule.
Each primary field is represented
by a Young diagram and we denote a primary field by the
corresponding Young diagram.

\begin{lemma}
\label{ff1}
We have the following fusion rule in the WZW-model
$SU(3)_{3k}$.
$$N_{ff}^{\unitlength 0.1mm
\begin{picture}(30,10)
\put(0,0){\line(1,0){30}}
\put(0,10){\line(1,0){30}}
\put(0,0){\line(0,1){10}}
\put(10,0){\line(0,1){10}}
\put(20,0){\line(0,1){10}}
\put(30,0){\line(0,1){10}}
\end{picture}}=1.$$
\end{lemma}

\begin{proof}
By \cite{GN}, we can apply the fusion rule described in
\cite{GW}.  By the Young--Pieri rule in
\cite[Proposition 2.6 (a)]{GW} and the classical
Littlewood--Richardson rule (see \cite[Section 1.9]{Mac},
for example), we get the conclusion.
\qed\end{proof}

\begin{lemma}
\label{ff2}
We have the following fusion rule in the WZW-model
$SU(3)_{3k}$.
$$N_{ff}^{\unitlength 0.1mm
\begin{picture}(20,20)
\put(0,0){\line(1,0){10}}
\put(0,10){\line(1,0){20}}
\put(0,20){\line(1,0){20}}
\put(0,0){\line(0,1){20}}
\put(10,0){\line(0,1){20}}
\put(20,10){\line(0,1){10}}
\end{picture}}=2.$$
\end{lemma}

\begin{proof}
This follows from Lemma \ref{ff1} and
$${\unitlength 0.2mm
\begin{picture}(11,10)
\put(0,0){\line(1,0){10}}
\put(0,10){\line(1,0){10}}
\put(0,0){\line(0,1){10}}
\put(10,0){\line(0,1){10}}
\end{picture}}^3
=\emptyset+\;
\unitlength 0.2mm
\begin{picture}(30,10)
\put(0,0){\line(1,0){30}}
\put(0,10){\line(1,0){30}}
\put(0,0){\line(0,1){10}}
\put(10,0){\line(0,1){10}}
\put(20,0){\line(0,1){10}}
\put(30,0){\line(0,1){10}}
\end{picture}
\;+2\;\;
\unitlength 0.2mm
\begin{picture}(20,20)
\put(0,0){\line(1,0){10}}
\put(0,10){\line(1,0){20}}
\put(0,20){\line(1,0){20}}
\put(0,0){\line(0,1){20}}
\put(10,0){\line(0,1){20}}
\put(20,10){\line(0,1){10}}
\end{picture},$$
because $f\;{\unitlength 0.2mm
\begin{picture}(11,10)
\put(0,0){\line(1,0){10}}
\put(0,10){\line(1,0){10}}
\put(0,0){\line(0,1){10}}
\put(10,0){\line(0,1){10}}
\end{picture}}^3$ contains 6 copies of $f$.
\qed\end{proof}

\begin{lemma}
\label{fr1}
We have the following  identity in the WZW-model
$SU(3)_{3k}$.
$$\sum_{\gr(a)=0} N_{a
{\unitlength 0.1mm
\begin{picture}(11,10)
\put(0,0){\line(1,0){10}}
\put(0,10){\line(1,0){10}}
\put(0,0){\line(0,1){10}}
\put(10,0){\line(0,1){10}}
\end{picture}}^3}^a=
\sum_{\gr(b)=1} N_{b
{\unitlength 0.1mm
\begin{picture}(11,10)
\put(0,0){\line(1,0){10}}
\put(0,10){\line(1,0){10}}
\put(0,0){\line(0,1){10}}
\put(10,0){\line(0,1){10}}
\end{picture}}^3}^b.$$
\end{lemma}

\begin{proof}
Since the level is $3k$, the numbers of the primary
fields of grade 0, 1, 2 are $3k(k+1)/2+1$, 
$3k(k+1)/2$, $3k(k+1)/2$ respectively.
Recall that the primary fields are arrayed in a triangular
picture for $SU(3)_k$ as in Figure \ref{proj-id}.

We have three primary fields of
grade 0 at the three corners of the triangle.  The contribution
of these terms on the left hand side of the identity in this
Lemma is 3.  (See \cite{GW}, \cite[Section 1.9]{Mac},
\cite{Wal} again for the computations of the fusion rule.)

We have $3(k-1)$ primary fields of grade 0 on the three
edges of the triangle with three corners excluded.
Each term gives a contribution of
3 on the left hand side of the identity,
so we get $9(k-1)$ as the total contribution.

We next have $3k^2/2-3k/2+1$ primary fields of grade 0
inside the triangle.  Each term gives a contribution of
6 on the left hand side of the identity,
so we get $9k^2-9k+6$ as the total contribution.

The sum of these three contributions is $9k^2$ and this
is the number on the left hand side of the identity.

We similarly evaluate the right hand side of the identity.

We have $3k$ primary fields of grade 1 on the three
edges of the triangle.  Each term gives a contribution of
3 on the right hand side of the identity,
so we get $9k$ as the total contribution.

We next have $3k^2/2-3k/2$ primary fields of grade 1
inside the triangle.  Each term gives a contribution of
6 on the right hand side of the identity,
so we get $9k^2-9k$ as the total contribution.

The sum of these two contributions is $9k^2$, which
is equal to the left hand side.
\qed\end{proof}

\begin{lemma}
\label{fr2}
We have the following  identity in the WZW-model
$SU(3)_{3k}$.
$$\sum_{\gr(a)=0} N_{a \bar a}^
{\unitlength 0.1mm
\begin{picture}(30,10)
\put(0,0){\line(1,0){30}}
\put(0,10){\line(1,0){30}}
\put(0,0){\line(0,1){10}}
\put(10,0){\line(0,1){10}}
\put(20,0){\line(0,1){10}}
\put(30,0){\line(0,1){10}}
\end{picture}}=
\sum_{\gr(b)=1} N_{b \bar b}^
{\unitlength 0.1mm
\begin{picture}(30,10)
\put(0,0){\line(1,0){30}}
\put(0,10){\line(1,0){30}}
\put(0,0){\line(0,1){10}}
\put(10,0){\line(0,1){10}}
\put(20,0){\line(0,1){10}}
\put(30,0){\line(0,1){10}}
\end{picture}}+1.
$$
\end{lemma}

\begin{proof}
Let $\a$ be
$\unitlength 0.2mm
\begin{picture}(30,10)
\put(0,0){\line(1,0){30}}
\put(0,10){\line(1,0){30}}
\put(0,0){\line(0,1){10}}
\put(10,0){\line(0,1){10}}
\put(20,0){\line(0,1){10}}
\put(30,0){\line(0,1){10}}
\end{picture}$.
We next count the number of the paths from the
odd vertex $0$ to the odd vertex labeled with
$\a$ on both the principal and the dual principal
graphs.
The number $\sum_{\gr(a)=0} N_{a \bar a}^\a$ gives the number
of the paths on 
the principal graph.

Let $l$ be the number of the edges connecting the
odd vertex $\a$ and the even vertex labeled with
$p_{f,f}^{(0)}$ on the dual principal graph.
By Lemmas \ref{proj-id11} and \ref{ff1}, we know that
$l$ is 0 or 1.  Lemmas \ref{proj-id10}, \ref{proj-id11},
and \ref{ff1} imply that the number of the paths connecting
$0$ and $\alpha$ on the dual principal graph is
$(\sum_b N_{b \bar b}^\a-1)/3+l$, where the
summation is over all the primary fields $b$ in the
model $SU(3)_{3k}$.

Since the two numbers of the paths are equal,
we get
$$2\sum_{\gr(a)=0} N_{a \bar a}^\a
=\sum_{\gr(b)=1,2} N_{b \bar b}^\a-1+3l=
2\sum_{\gr(b)=1} N_{b \bar b}^\a-1+3l.$$
This implies $l=1$ because the both sides are
even numbers.  We then get the conclusion.
\qed\end{proof}

\begin{lemma}
\label{fr3}
We have the following  identity in the WZW-model
$SU(3)_{3k}$.
$$\sum_{\gr(a)=0}
N_{a \bar a}^{\unitlength 0.1mm
\begin{picture}(20,20)
\put(0,0){\line(1,0){10}}
\put(0,10){\line(1,0){20}}
\put(0,20){\line(1,0){20}}
\put(0,0){\line(0,1){20}}
\put(10,0){\line(0,1){20}}
\put(20,10){\line(0,1){10}}
\end{picture}}=
\sum_{\gr(b)=1}
N_{b \bar b}^{\unitlength 0.1mm
\begin{picture}(20,20)
\put(0,0){\line(1,0){10}}
\put(0,10){\line(1,0){20}}
\put(0,20){\line(1,0){20}}
\put(0,0){\line(0,1){20}}
\put(10,0){\line(0,1){20}}
\put(20,10){\line(0,1){10}}
\end{picture}}-1.$$
\end{lemma}

\begin{proof}
Recall that we have
\begin{equation}
\label{id1}
{\unitlength 0.2mm
\begin{picture}(11,10)
\put(0,0){\line(1,0){10}}
\put(0,10){\line(1,0){10}}
\put(0,0){\line(0,1){10}}
\put(10,0){\line(0,1){10}}
\end{picture}}^3
=\emptyset+\;
\unitlength 0.2mm
\begin{picture}(30,10)
\put(0,0){\line(1,0){30}}
\put(0,10){\line(1,0){30}}
\put(0,0){\line(0,1){10}}
\put(10,0){\line(0,1){10}}
\put(20,0){\line(0,1){10}}
\put(30,0){\line(0,1){10}}
\end{picture}
\;+2\;\;
\unitlength 0.2mm
\begin{picture}(20,20)
\put(0,0){\line(1,0){10}}
\put(0,10){\line(1,0){20}}
\put(0,20){\line(1,0){20}}
\put(0,0){\line(0,1){20}}
\put(10,0){\line(0,1){20}}
\put(20,10){\line(0,1){10}}
\end{picture}.
\end{equation}

Lemma \ref{fr1} and Frobenius reciprocity imply
\begin{equation}
\label{id2}
\sum_{\gr(a)=0} N_{a \bar a}^
{{\unitlength 0.1mm
\begin{picture}(11,10)
\put(0,0){\line(1,0){10}}
\put(0,10){\line(1,0){10}}
\put(0,0){\line(0,1){10}}
\put(10,0){\line(0,1){10}}
\end{picture}}^3}=
\sum_{\gr(b)=1} N_{b \bar b}^
{{\unitlength 0.1mm
\begin{picture}(11,10)
\put(0,0){\line(1,0){10}}
\put(0,10){\line(1,0){10}}
\put(0,0){\line(0,1){10}}
\put(10,0){\line(0,1){10}}
\end{picture}}^3}.
\end{equation}

We also have the easy identity,
\begin{equation}
\label{id3}
\sum_{\gr(a)=0} N_{a \bar a}^\emptyset=
\sum_{\gr(b)=0} N_{b \bar b}^\emptyset+1,
\end{equation}
since the both  sides are equal to $3k(k+1)/2+1$.

Identities (\ref{id1}), (\ref{id2}), (\ref{id3}) and
Lemma \ref{fr2} imply the conclusion.
\qed\end{proof}

\begin{lemma}
\label{edge0}
The odd vertex of the dual principal graph labeled with
$\unitlength 0.2mm
\begin{picture}(20,20)
\put(0,0){\line(1,0){10}}
\put(0,10){\line(1,0){20}}
\put(0,20){\line(1,0){20}}
\put(0,0){\line(0,1){20}}
\put(10,0){\line(0,1){20}}
\put(20,10){\line(0,1){10}}
\end{picture}$
is not connected to the
even vertex labeled with $p^{(0)}_{f,f}$.
\end{lemma}

\begin{proof}
We count the number of appropriate paths again.

The number of paths connecting the odd vertex
$0$ to
$\unitlength 0.2mm
\begin{picture}(20,20)
\put(0,0){\line(1,0){10}}
\put(0,10){\line(1,0){20}}
\put(0,20){\line(1,0){20}}
\put(0,0){\line(0,1){20}}
\put(10,0){\line(0,1){20}}
\put(20,10){\line(0,1){10}}
\end{picture}$
on the principal graph of the asymptotic
inclusion is
$\sum_{\gr(a)=0}
N_{a \bar a}^{\unitlength 0.1mm
\begin{picture}(20,20)
\put(0,0){\line(1,0){10}}
\put(0,10){\line(1,0){20}}
\put(0,20){\line(1,0){20}}
\put(0,0){\line(0,1){20}}
\put(10,0){\line(0,1){20}}
\put(20,10){\line(0,1){10}}
\end{picture}}$
because the principal graph is the fusion graph.

Let $l$ be the number of edges connecting
the odd vertex of the dual principal graph labeled with
$\unitlength 0.2mm
\begin{picture}(20,20)
\put(0,0){\line(1,0){10}}
\put(0,10){\line(1,0){20}}
\put(0,20){\line(1,0){20}}
\put(0,0){\line(0,1){20}}
\put(10,0){\line(0,1){20}}
\put(20,10){\line(0,1){10}}
\end{picture}$
to the even vertex labeled with $p^{(0)}_{f,f}$.

The number of paths connecting the odd vertex
$0$ to
$\unitlength 0.2mm
\begin{picture}(20,20)
\put(0,0){\line(1,0){10}}
\put(0,10){\line(1,0){20}}
\put(0,20){\line(1,0){20}}
\put(0,0){\line(0,1){20}}
\put(10,0){\line(0,1){20}}
\put(20,10){\line(0,1){10}}
\end{picture}$
on the dual principal graph of the asymptotic
inclusion is $\left(
\sum_{b}
N_{b \bar b}^{\unitlength 0.1mm
\begin{picture}(20,20)
\put(0,0){\line(1,0){10}}
\put(0,10){\line(1,0){20}}
\put(0,20){\line(1,0){20}}
\put(0,0){\line(0,1){20}}
\put(10,0){\line(0,1){20}}
\put(20,10){\line(0,1){10}}
\end{picture}}
-
N_{ff}^{\unitlength 0.1mm
\begin{picture}(20,20)
\put(0,0){\line(1,0){10}}
\put(0,10){\line(1,0){20}}
\put(0,20){\line(1,0){20}}
\put(0,0){\line(0,1){20}}
\put(10,0){\line(0,1){20}}
\put(20,10){\line(0,1){10}}
\end{picture}}
\right)/3+l$.
Lemmas \ref{ff2}, \ref{fr3} show $l=0$.
\qed\end{proof}

We finally prove the main theorem in this Section as follows.

\begin{theorem}
For the subfactor $N\subset M$ arising from the WZW-model
$SU(3)_{3k}$, 
Assumption \ref{orbif} holds.
\end{theorem}

\begin{proof}
Since we have Lemmas \ref{proj-id11}, \ref{ff2},
we have one of the following two cases.
\begin{enumerate}
\item We have a minimal central projection $p_{ff}^{(1)}$
majorized by $p_{ff}/3$ such that the odd vertex labeled with
$\unitlength 0.2mm
\begin{picture}(20,20)
\put(0,0){\line(1,0){10}}
\put(0,10){\line(1,0){20}}
\put(0,20){\line(1,0){20}}
\put(0,0){\line(0,1){20}}
\put(10,0){\line(0,1){20}}
\put(20,10){\line(0,1){10}}
\end{picture}$ is connected to the even vertex labeled with
$p_{ff}^{(1)}$ on the dual principal graph of the asymptotic
inclusion by exactly two edges.
\item We have two minimal central projections
$p_{ff}^{(1)}$, $p_{ff}^{(2)}$
majorized by $p_{ff}/3$ such that the odd vertex labeled with
$\unitlength 0.2mm
\begin{picture}(20,20)
\put(0,0){\line(1,0){10}}
\put(0,10){\line(1,0){20}}
\put(0,20){\line(1,0){20}}
\put(0,0){\line(0,1){20}}
\put(10,0){\line(0,1){20}}
\put(20,10){\line(0,1){10}}
\end{picture}$ is connected to each of
the even vertices labeled with
$p_{ff}^{(1)}$,  $p_{ff}^{(2)}$
on the dual principal graph of the asymptotic
inclusion by exactly one edge.
\end{enumerate}
Suppose that we have Case 1.  Lemma \ref{PF}
shows that the dimension of the bimodule corresponding
to $p_{ff}^{(1)}$ is equal to that of the bimodule
corresponding to $p_{ff}^{(0)}$.  Lemma \ref{glob}
implies that the central projection
$p_{ff}^{(2)}=p_{ff}/3-p_{ff}^{(0)}-p_{ff}^{(1)}$
is minimal and the dimension of the bimodule corresponding
to $p_{ff}^{(2)}$ is also equal to that of the bimodule
corresponding to $p_{ff}^{(0)}$.

Next suppose that we have Case 2.  Lemma \ref{PF}
implies that the sum of the dimensions of the
bimodules  corresponding
to $p_{ff}^{(1)}$, $p_{ff}^{(2)}$ is equal to twice
of that of the bimodule
corresponding to $p_{ff}^{(0)}$.
This shows
$p_{ff}/3=p_{ff}^{(0)}+p_{ff}^{(1)}+p_{ff}^{(2)}$.
Then Lemma \ref{glob} then
implies that these two dimensions have to be equal.

In any case, we have a decomposition 
$p_{ff}/3=p_{ff}^{(0)}+p_{ff}^{(1)}+p_{ff}^{(2)}$ into
minimal central projections and each of the three minimal
central projection has the same corresponding dimension.
This completes the proof.
\qed\end{proof}

This Theorem implies the following by a simple computation.
This Corollary is a generalized version of Corollary \ref{ocn2}.
Again note that this Corollary
gives the number of even vertices of the dual principal
graph of the asymptotic inclusions and that these are also the dimensions
of the Hilbert spaces $H_{S^1\times S^1}$ in the corresponding
topological quantum field theories for the original subfactors.

\begin{corollary}
Let $N\subset M$ be the subfactor corresponding to
$SU(3)_k$ with $k>2$.  Then the number of
the irreducible $M_\infty$-$M_\infty$ bimodules arising
from the asymptotic inclusion is given as follows.
$$\begin{array}{ll}
\displaystyle\frac{(k+1)^2(k+2)^2}{36},
&\quad {\rm if \ }k\not\equiv 0{\ mod\ }3,\\
\displaystyle\frac{k^4+6k^3+13k^2+12k+108}{36},
&\quad {\rm if \ }k\equiv 0{\ mod\ }3.
\end{array}$$
\end{corollary}

As examples, we work out the dual principal graphs 
for small $k$ such as $k=3,6$ in the rest of this Section.

First, we label the primary fields of $SU(3)_3$ as in
Figure \ref{su33}.

\unitlength 0.2mm
\thinlines
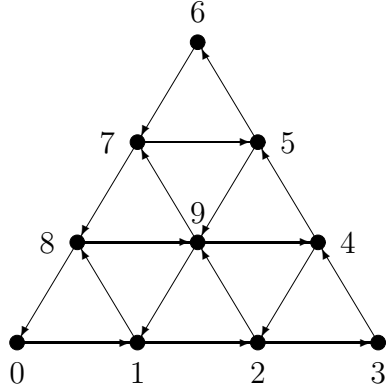
\begin{figure}[tb]
\begin{center}
\begin{picture}(300,270)(0,0)
\multiput(105,30)(80,0){3}{\vector(1,0){0}}
\multiput(145,96.67)(80,0){2}{\vector(1,0){0}}
\put(185,163.33){\vector(1,0){0}}
\multiput(32.5,34.17)(80,0){3}{\vector(-1,-2){0}}
\multiput(72.5,100.84)(80,0){2}{\vector(-1,-2){0}}
\put(112.5,167.5){\vector(-1,-2){0}}
\multiput(72.5,92.5)(80,0){3}{\vector(-1,2){0}}
\multiput(112.5,159.17)(80,0){2}{\vector(-1,2){0}}
\put(152.5,225.83){\vector(-1,2){0}}
\put(30,30){\line(1,0){80}}
\put(110,30){\line(1,0){80}}
\put(190,30){\line(1,0){80}}
\put(70,96.67){\line(1,0){80}}
\put(150,96.67){\line(1,0){80}}
\put(110,163.33){\line(1,0){80}}
\put(110,30){\line(-3,5){40}}
\put(190,30){\line(-3,5){40}}
\put(270,30){\line(-3,5){40}}
\put(150,96.67){\line(-3,5){40}}
\put(230,96.67){\line(-3,5){40}}
\put(190,163.33){\line(-3,5){40}}
\put(70,96.67){\line(-3,-5){40}}
\put(150,96.67){\line(-3,-5){40}}
\put(230,96.67){\line(-3,-5){40}}
\put(110,163.33){\line(-3,-5){40}}
\put(190,163.33){\line(-3,-5){40}}
\put(150,230){\line(-3,-5){40}}
\put(30,30){\circle*{10}}
\put(110,30){\circle*{10}}
\put(190,30){\circle*{10}}
\put(270,30){\circle*{10}}
\put(70,96.67){\circle*{10}}
\put(150,96.67){\circle*{10}}
\put(230,96.67){\circle*{10}}
\put(110,163.33){\circle*{10}}
\put(190,163.33){\circle*{10}}
\put(150,230){\circle*{10}}
\put(30,10){\makebox(0,0){$0$}}
\put(110,10){\makebox(0,0){$1$}}
\put(190,10){\makebox(0,0){$2$}}
\put(270,10){\makebox(0,0){$3$}}
\put(50,96.67){\makebox(0,0){$8$}}
\put(150,116.67){\makebox(0,0){$9$}}
\put(250,96.67){\makebox(0,0){$4$}}
\put(90,163.33){\makebox(0,0){$7$}}
\put(210,163.33){\makebox(0,0){$5$}}
\put(150,250){\makebox(0,0){$6$}}
\end{picture}
\end{center}
\caption{Primary fields for $SU(3)_3$}
\label{su33}
\end{figure}

Then the principal graph of the asymptotic inclusion
of the subfactor corresponding to $SU(3)_3$ is given
as the fusion graph as in the upper half of Figure
\ref{su33-2}.  For the dual principal graph,
we know the graph except for the edges connected to the
three vertices $(99)_0, (99)_1, (99)_2$.  From the
Perron--Frobenius property, we can determine these edges
as in the bottom half of Figure \ref{su33-2}.  These
edges are marked thick.

Since the subfactor corresponding to $SU(3)_3$ has
index 4 and is described as $R\rtimes A_3\subset R\rtimes A_4$,
where $A_3$ and $A_4$ are the alternating groups of order
3 and 4 respectively and these groups act freely on
the hyperfinite II$_1$ factor $R$,  the paragroup
of the asymptotic inclusion is given by that of the
subfactor $R^{A_4\times A_4}\subset R^{A_4}$, where
$A_4$ is diagonally embedded into $A_4\times A_4$ and
the group $A_4$ acts freely on $R$, by Ocneanu's theorem
again.  (See \cite[Lemma 2.15]{K2}, \cite[Appendix]{K3},
\cite[Section 12.8]{EK4}.)

So the (dual) principal graphs of the asymptotic inclusion
can be described with Ocneanu's theorem again.  (See \cite{KY}.)
Of course, this method gives the same result as in
Figure \ref{su33-2}.

\unitlength 0.47mm
\thinlines
\begin{figure}[tb]
\begin{center}
\begin{picture}(320,150)
\multiput(10,140)(20,0){16}{\circle*{2}}
\multiput(10,20)(20,0){14}{\circle*{2}}
\put(30,80){\circle*{2}}
\put(90,80){\circle*{2}}
\put(150,80){\circle*{2}}
\put(230,80){\circle*{2}}
\path(30,80)(10,140)
\path(30,80)(30,140)
\path(30,80)(50,140)
\path(90,80)(70,140)
\path(90,80)(90,140)
\path(90,80)(110,140)
\path(150,80)(130,140)
\path(150,80)(150,140)
\path(150,80)(170,140)
\path(150,80)(190,140)
\path(30,80)(190,140)
\path(90,80)(190,140)
\path(231,80)(191,140)
\path(229,80)(189,140)
\path(230,80)(210,140)
\path(230,80)(230,140)
\path(230,80)(250,140)
\path(230,80)(270,140)
\path(230,80)(290,140)
\path(230,80)(310,140)
\path(30,80)(10,20)
\path(30,80)(30,20)
\path(30,80)(50,20)
\path(90,80)(70,20)
\path(90,80)(130,20)
\path(90,80)(150,20)
\path(150,80)(110,20)
\path(150,80)(170,20)
\path(150,80)(190,20)
\path(230,80)(30,20)
\path(230,80)(50,20)
\path(230,80)(130,20)
\path(230,80)(150,20)
\path(230,80)(170,20)
\path(230,80)(190,20)
\path(230,80)(210,20)
\path(230,80)(230,20)
\Thicklines
\path(30,80)(90,20)
\path(90,80)(90,20)
\path(150,80)(90,20)
\path(230,80)(250,20)
\path(230,80)(270,20)
\put(3,15){\makebox(0,0){$*$}}
\put(10,15){\makebox(0,0){$00$}}
\put(10,8){\makebox(0,0){$36$}}
\put(10,1){\makebox(0,0){$63$}}
\put(30,15){\makebox(0,0){$81$}}
\put(30,8){\makebox(0,0){$54$}}
\put(30,1){\makebox(0,0){$27$}}
\put(50,15){\makebox(0,0){$18$}}
\put(50,8){\makebox(0,0){$45$}}
\put(50,1){\makebox(0,0){$72$}}
\put(70,15){\makebox(0,0){$03$}}
\put(70,8){\makebox(0,0){$30$}}
\put(70,1){\makebox(0,0){$66$}}
\put(90,15){\makebox(0,0){$(99)_0$}}
\put(110,15){\makebox(0,0){$06$}}
\put(110,8){\makebox(0,0){$33$}}
\put(110,1){\makebox(0,0){$60$}}
\put(130,15){\makebox(0,0){$21$}}
\put(130,8){\makebox(0,0){$84$}}
\put(130,1){\makebox(0,0){$57$}}
\put(150,15){\makebox(0,0){$12$}}
\put(150,8){\makebox(0,0){$48$}}
\put(150,1){\makebox(0,0){$75$}}
\put(170,15){\makebox(0,0){$51$}}
\put(170,8){\makebox(0,0){$24$}}
\put(170,1){\makebox(0,0){$87$}}
\put(190,15){\makebox(0,0){$15$}}
\put(190,8){\makebox(0,0){$42$}}
\put(190,1){\makebox(0,0){$78$}}
\put(210,15){\makebox(0,0){$09$}}
\put(210,8){\makebox(0,0){$39$}}
\put(210,1){\makebox(0,0){$69$}}
\put(230,15){\makebox(0,0){$90$}}
\put(230,8){\makebox(0,0){$93$}}
\put(230,1){\makebox(0,0){$96$}}
\put(250,15){\makebox(0,0){$(99)_1$}}
\put(270,15){\makebox(0,0){$(99)_2$}}
\put(23,80){\makebox(0,0){$0$}}
\put(83,80){\makebox(0,0){$3$}}
\put(143,80){\makebox(0,0){$6$}}
\put(220,80){\makebox(0,0){$9$}}
\put(3,145){\makebox(0,0){$*$}}
\put(10,145){\makebox(0,0){$00$}}
\put(30,145){\makebox(0,0){$36$}}
\put(50,145){\makebox(0,0){$63$}}
\put(70,145){\makebox(0,0){$03$}}
\put(90,145){\makebox(0,0){$30$}}
\put(110,145){\makebox(0,0){$66$}}
\put(130,145){\makebox(0,0){$06$}}
\put(150,145){\makebox(0,0){$33$}}
\put(170,145){\makebox(0,0){$60$}}
\put(190,145){\makebox(0,0){$99$}}
\put(210,145){\makebox(0,0){$09$}}
\put(230,145){\makebox(0,0){$39$}}
\put(250,145){\makebox(0,0){$69$}}
\put(270,145){\makebox(0,0){$90$}}
\put(290,145){\makebox(0,0){$93$}}
\put(310,145){\makebox(0,0){$96$}}
\end{picture}
\end{center}
\caption{(dual) principal graphs for $SU(3)_3$}
\label{su33-2}
\end{figure}
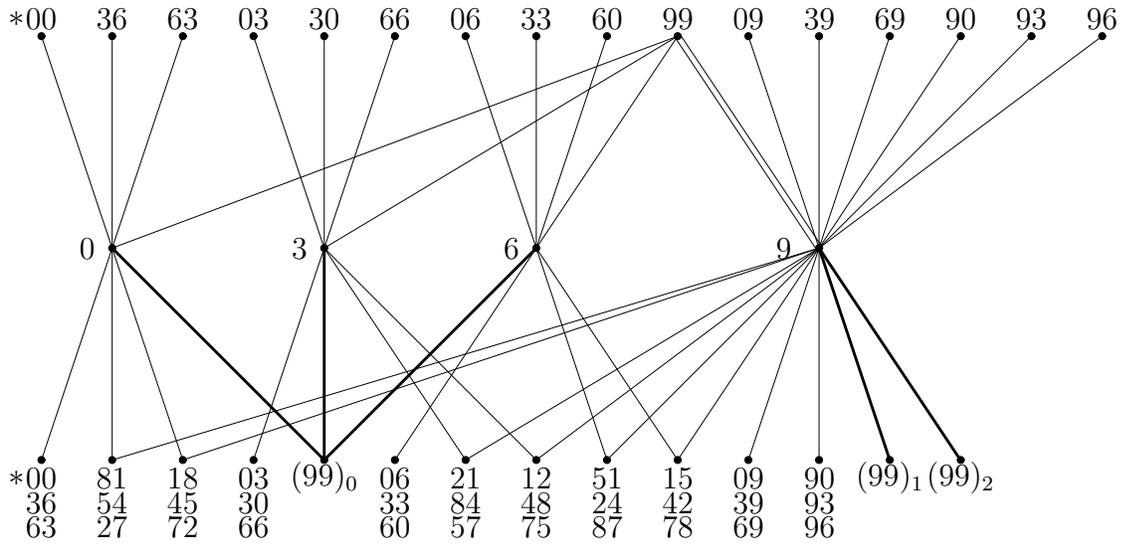

The next example is $SU(3)_6$.  In this case, the system
${\cal M}$ has 10 primary fields
and thus the principal graph of the
asymptotic inclusion has 100 even vertices, and the
dual principal graph has 90 even vertices.  Since these
graphs are too complicated, we draw only the edges concerned
with the three even vertices $p_{ff}^{(0)}$, 
$p_{ff}^{(1)}$, $p_{ff}^{(2)}$.  Then the Perron--Frobenius
property and counting of paths with unitarity gives the graph
as in Figure \ref{su36}.  In this Figure, the symbol
$(lm)$ denotes the Young diagram with $l$ boxes in the first
row and $m$ boxes in the second row.

\unitlength 0.75mm
\thinlines
\begin{figure}[tb]
\begin{center}
\begin{picture}(200,50)
\multiput(10,40)(20,0){10}{\circle*{2}}
\put(70,10){\circle*{2}}
\put(150,10){\circle*{2}}
\put(190,10){\circle*{2}}
\path(70,10)(10,40)
\path(70,10)(30,40)
\path(70,10)(50,40)
\path(70,10)(70,40)
\path(70,10)(90,40)
\path(70,10)(110,40)
\path(70,10)(130,40)
\path(150,10)(130,40)
\path(150,10)(150,40)
\path(150,10)(170,40)
\path(150,10)(190,40)
\path(190,10)(130,40)
\path(190,10)(150,40)
\path(190,10)(170,40)
\path(190,10)(190,40)
\put(70,5){\makebox(0,0){$(42)(42)_0$}}
\put(150,5){\makebox(0,0){$(42)(42)_1$}}
\put(190,5){\makebox(0,0){$(42)(42)_2$}}
\put(10,45){\makebox(0,0){$(0)$}}
\put(30,45){\makebox(0,0){$(6)$}}
\put(50,45){\makebox(0,0){$(66)$}}
\put(70,45){\makebox(0,0){$(3)$}}
\put(90,45){\makebox(0,0){$(33)$}}
\put(110,45){\makebox(0,0){$(63)$}}
\put(130,45){\makebox(0,0){$(42)$}}
\put(150,45){\makebox(0,0){$(21)$}}
\put(170,45){\makebox(0,0){$(51)$}}
\put(190,45){\makebox(0,0){$(54)$}}
\end{picture}
\end{center}
\caption{Part of the dual principal graphs for $SU(3)_6$}
\label{su36}
\end{figure}
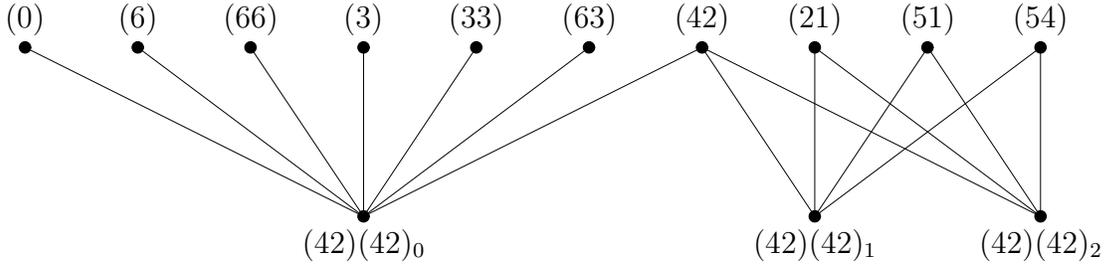

\section{Orbifold subfactors}
\label{sect-orbif}

In Sections \ref{sect-su2}, \ref{sect-su3}, we have observed that
the even vertices of the dual principal graphs of the
asymptotic inclusions are given by merging/splitting of the
vertices with symmetries on pairs of the original labels.
In the $SU(2)_k$ case, Ocneanu
has noticed that this situation is similar to the orbifold
construction for subfactors studied by us in \cite{EK1}, \cite{K1}.
(See also \cite{G}, \cite{X}.)
However, the dual principal graphs we have studied in
Sections \ref{sect-su2}, \ref{sect-su3} are {\sl not} orbifold
graphs in the sense of \cite{EK1}, \cite{K1}, \cite{X},
because we have merging/splitting of the
vertices only for the even vertices.  In this Section, we study
a relation of this orbifold phenomena of Ocneanu to the
orbifold construction in our sense.

Let $N\subset M$ be the Jones subfactor of type $A_{4n-3}$.
That is, it is the hyperfinite type II$_1$ subfactor
corresponding to $SU(2)_{4n-4}$.  To avoid disconnectedness
of the fusion graph, we assume that $n>2$.  (If $n=2$, we get
a subfactor arising from a free action of a group $S_3$, so
everything can be studied with classical methods on group
actions.)  As in \cite{K1},
we get the orbifold subfactor
$P=N\rtimes_\sigma \Z/2\Z\subset Q=M\rtimes_\sigma \Z/2\Z$
of type $D_{2n}$, where $\sigma$ gives a non-strongly-outer
action of $\Z/2\Z$ on the subfactor $N\subset M$ in the
sense of \cite{CK}.  (See also \cite{G}.)

Let $\a$ be the dual action of $\sigma$ on $P\subset Q$.  Then
we have $N=P^\a$ and $M=Q^\a$, of course.
Then the asymptotic inclusion
$M\vee (M'\cap M_\infty)\subset M_\infty$ is described as
$Q^\a\vee (Q'\cap Q_\infty)^\a\subset Q_\infty^\a$.
Putting $R=(Q\vee (Q'\cap Q_\infty))^\a$, we get
$Q^\a\vee (Q'\cap Q_\infty)^\a\subset R\subset Q_\infty^\a$ and
$[R: Q^\a\vee (Q'\cap Q_\infty)^\a]=2$.  This intermediate
subfactor corresponds to the
intermediate subfactor $(M_\omega)^\sigma$
of the central sequence subfactor
$N^\omega\cap M'\subset M_\omega$ described in
\cite[Section 3]{K2}, \cite[Section 4]{K3} in the
correspondence of Ocneanu \cite[page 42]{O2}, \cite[Theorem 4.1]{K3}.
(Here $\omega$ is a free ultrafilter over $\N$.
See also \cite[Theorem 15.32]{EK4}.)

We use the notation $[[M:N]]$ for the global index of $N\subset M$
as in \cite{S1}.
We easily get
$$[[M_\infty: M\vee (M'\cap M_\infty)]]/4
=[[Q_\infty: Q\vee (Q'\cap Q_\infty)]]$$
from the description of the
principal graph as the fusion graph.
Note that $R\subset Q_\infty^\a$ is given as
the simultaneous fixed point
algebras of $Q\vee (Q'\cap Q_\infty)\subset Q_\infty$ by the action
$\a$.  By looking at $\a$, we can conclude that
we have one of the following three cases.
\begin{enumerate}
\item $[[M_\infty:R]]=[[M_\infty: M\vee (M'\cap M_\infty)]]/4$.
\item $[[M_\infty:R]]=[[M_\infty: M\vee (M'\cap M_\infty)]]/2$.
\item $[[M_\infty:R]]=[[M_\infty: M\vee (M'\cap M_\infty)]]/8$.
\end{enumerate}
That is, if the action is strongly outer in the sense of
\cite{CK} and has a trivial Loi invariant, then we get Case 1,
if the action is strongly outer and has a non-trivial Loi
invariant, then we have Case 2,
and if the action is not strongly outer, then we have Case 3.
Since the fusion rule algebra of the $M_\infty$-$M_\infty$
bimodules arising from $R\subset M_\infty$ is a fusion rule
subalgebra of those arising from $M\vee(M'\cap M_\infty)
\subset M_\infty$ (see \cite[Lemma 2.4]{S1}, for example), 
we look for a fusion rule subalgebra of that of the 
$M_\infty$-$M_\infty$ bimodules arising from the asymptotic
inclusion of $N\subset M$.

We first study fusion rule subalgebras of the WZW-model
$SU(2)_{2k}$.

\begin{lemma}
\label{fusion-sub}
Let $\cal N$ be a closed subsystem of primary fields under fusion
of the WZW-model $SU(2)_{2k}$ labeled as $\{0,1,2,\dots,2k\}$. 
Then ${\cal N}$ is one of the following;
$\{0\}$, $\{0,2k\}$, $\{0,2,4,\dots,2k\}$, $\{0,1,2,\dots,2k\}$
\end{lemma}

\begin{proof}
It is clear that these four indeed give subsystems.

Suppose that ${\cal N}\neq\{0\}$.
Let $l$ be the smallest non-zero label appearing in $\cal N$.
If $l=1$, $l=2$, or $l=2k$, then we clearly
have ${\cal N}=\{0,1,2,\dots,2k\}$,
${\cal N}=\{0,2,4,\dots,2k\}$, ${\cal N}=\{0, 2k\}$, respectively.
If $2<l<2k$, then we would have $N_{ll}^2=1$, which implies $2
\in{\cal N}$ and thus a contradiction.
\qed\end{proof}

\begin{lemma}
\label{correct-gb}
Let ${\cal N}$ be the system of $M_\infty$-$M_\infty$
bimodules arising from the asymptotic inclusion of the subfactor
$N\subset M$ of type $A_{4n-3}$.
Let $\gamma$ be the global index of this system.
Suppose we have a subsystem ${\cal N}_0$ of $\cal N$ with global index
equal to one of $\gamma/2$, $\gamma/4$, $\gamma/8$.  Then
${\cal N}_0$ is a subsystem of the
$M_\infty$-$M_\infty$ bimodules labeled with pairs of
even numbers as in Section \ref{sect-su2} and its global index
is $\gamma/2$.
\end{lemma}

\begin{proof}
It is clear that the subsystem of the
$M_\infty$-$M_\infty$ bimodules labeled with pairs of
even numbers has global index $\gamma/2$.

Suppose that ${\cal N}_0$ contains a bimodule labeled with
a pair of odd numbers.  By taking an appropriate tensor
power of this bimodule, we have $(2,2)$ in this system ${\cal N}_0$.

We set $X$ be the set of
labels of pairs of integers appearing in ${\cal N}_0$ and
set $Y=\{l\mid (0,l)\in X\}$.  Then $Y$ gives a subsystem
of the original WZW-model $SU(2)_{4n-4}$.
By Lemma \ref{fusion-sub}, we have four cases for $Y$.
The assumption on the global index forces
$Y=\{0,2,4,\dots,4n-4\}$ and $[{\cal N}_0]=\gamma/2$.
Then we have the conclusion.
\qed\end{proof}

Let $\be$ be the $M_\infty$-$M_\infty$ bimodule labeled
with $(0,4n-4)=(4n-4,0)$.  It is clear that this bimodule
has dimension 1.
We can apply the orbifold construction for
tensor categories as in \cite{Y} and get the following Lemma.

\begin{lemma}
\label{orbif-d}
Let $N\subset M$, $P\subset Q$ be as above.
Let ${\cal N}_0$ be the system of $M_\infty$-$M_\infty$
bimodules arising from $R\subset Q_\infty^\a=M_\infty$.
Let ${\cal N}_1$ be the system of $Q_\infty$-$Q_\infty$ 
bimodules arising from the asymptotic inclusion of the subfactor
$P\subset Q$ of type $D_{2n}$.
Then the system ${\cal N}_1$ is given as the orbifold construction
of ${\cal N}_0$ with $\be$ as above.
\end{lemma}

\begin{proof}
Let $\cal N$ be the system of $M_\infty$-$M_\infty$
bimodules arising from the asymptotic inclusion of the subfactor
$N\subset M$ of type $A_{4n-3}$.

Lemma \ref{correct-gb} implies that $[{\cal N}_1]=[{\cal N}_0]/2$,
which gives the conclusion.
\qed\end{proof}

\begin{theorem}
\label{EKO}
Let $N\subset M$ be the Jones subfactor of type $A_{4n-3}$ with
$n>2$.
Let ${\cal N}_0$ be the subsystem of $M_\infty$-$M_\infty$ bimodules
arising from the asymptotic inclusion $M\vee (M'\cap M_\infty)
\subset M_\infty$ labeled with pairs of even integers as in Section
\ref{sect-su2}.

Let $\sigma$ be the outer, non-strongly-outer automorphism of order 2
of $N\subset M$.  The system ${\cal N}_0$ is isomorphic to the
system of
$(M\otimes M)^{\si\otimes\si}$-$(M\otimes M)^{\si\otimes\si}$
bimodules arising from the orbifold subfactor
$(N\otimes N)^{\si\otimes\si}\subset (M\otimes M)^{\si\otimes\si}$.
\end{theorem}

\begin{proof}
This follows from Lemma \ref{orbif-d}.
\qed\end{proof}

The meaning of the above Theorem is as follows.
When we apply the ``quantum double'' construction to
a degenerate system, it is not enough to take a simple
``double'' to because of degeneracy.  Pairs labeled with
ghosts appear so that the non-degeneracy is recovered,
but then we have too many bimodules and the global index,
giving the size of the system, becomes too large.  Then
the orbifold construction removes this redundancy and the
correct global index is realized.  The bimodules labeled
with pairs of ghosts disappear when we remove an intermediate
subfactor of index 2, which is the order of the orbifold
construction.

\section{Orbifold construction for braiding}
\label{sect-orbif-braid}

\begin{theorem}
\label{D2n-braid}
The system of the $M$-$M$ bimodules arising from a subfactor
$N\subset M$ of type $D_{2n}$, $n>2$, has a non-degenerate braiding.
\end{theorem}

\begin{proof}
The system of the $Q_\infty$-$Q_\infty$ bimodules has a non-degenerate
braiding by Ocneanu's general theory.
(See \cite[Section 12.7]{EK4}, for example.)

Lemma \ref{orbif-d} implies that the system given by the
orbifold construction on the system
$(0,0), (0,2),\dots, (0,4n-4)=\be$ with $\be$ is a subsystem of
the $Q_\infty$-$Q_\infty$ bimodules.  We thus get a braiding
naturally.  The non-degeneracy is also easy to see, because
if have degeneracy, then the degenerate subsystem would give
a finite abelian group by \cite{DR}, which is impossible
by $n>2$.
\qed\end{proof}

\begin{corollary}
\label{D2n-asymp}
The dual principal graph of the asymptotic inclusion of the
hyperfinite II$_1$ subfactor $N\subset M$
with principal graph $D_{2n}$  is
the fusion graph of the system of  $M$-$M$ bimodules.
\end{corollary}

\begin{proof}
This follows from Theorem \ref{D2n-braid} and
Proposition \ref{edges}.
\qed\end{proof}

\begin{remark}
\label{others}
{\rm Ocneanu has constructed a braiding on the even vertices of
$D_{2n}$ with an entirely different method in \cite{O8}.
His theory in \cite{O8} also shows that his
braiding and ours must be the same.

Turaev and Wenzl \cite{TW} have worked on a similar construction
to our orbifold construction in categories of tangles.
In their approach to the Reshetikhin--Turaev type topological
quantum field theory \cite{RT}, they need a certain non-degeneracy
and make some construction similar to our orbifold construction
to remove the degeneracy.  It
seems that their construction, in particular, gives a 
braiding on the even vertices of $D_{2n}$ and we expect that
their braiding is also same as ours, but the actual relation
is not clear.

The basic idea is that the orbifold construction can be performed
when we have some kind of degeneracy and this degeneracy is removed
by the orbifold construction.}
\end{remark}

\end{document}